\documentstyle[12pt]{article}
\newcommand{\sect}[1]{\section{#1}\setcounter{equation}{0}}

\font\mbn=msbm10 scaled \magstep1
\font\mbs=msbm7 scaled \magstep1
\font\mbss=msbm5 scaled \magstep1
\newfam\mbff
\textfont\mbff=\mbn
\scriptfont\mbff=\mbs
\scriptscriptfont\mbff=\mbss\def\mbf{\fam\mbff}
\def\Re{{\mbf R}}

\def\Co{{\mbf C}}

\def\H{{\mbf H}}
\newtheorem{Th}{Theorem}[section]
\newtheorem{Lm}[Th]{Lemma}
\newtheorem{C}[Th]{Corollary}

\newtheorem{D}[Th]{Definition}
\newtheorem{Proposition}[Th]{Proposition}
\newtheorem{R}[Th]{Remark}
\newtheorem{E}[Th]{Example}
\author{Alexander Brudnyi\thanks{Research supported in part by NSERC.
}\\
Department of Mathematics and Statistics\\
University of Calgary, Calgary\\
Canada\\
\\
Yuri Brudnyi\\
Department of Mathematics\\
Technion, Haifa\\
Israel}
\title{A Universal Lipschitz Extension Property of Gromov Hyperbolic Spaces}
\date{}
\begin{document}
\maketitle
\begin{abstract}
{A metric space has the universal Lipschitz extension property if for each  
subspace $S$ embedded
quasi-isometrically into an arbitrary metric space $M$ there exists a continuous linear
extension of Banach-valued Lipschitz functions on $S$ to those on all of
$M$. We show that the 
finite direct sum of Gromov hyperbolic spaces of bounded geometry is universal in the sense of this definition.}
\end{abstract}
\sect{Formulation of Main Results}
In order to present a precise formulation of the main results we need several definitions.

Let $(M,d)$ be a metric space with underlying set $M$ and metric $d$ (we write simply $M$ if
$d$ can be restored from the context). The space of Banach-valued Lipschitz functions on $M$
with target space $X$ is denoted by $Lip(M,X)$; this space is endowed with the standard seminorm
\begin{equation}\label{e1}
L(f):=\sup_{m\neq m'}\left\{\frac{||f(m)-f(m')||}{d(m,m')}\right\}.
\end{equation}
(In case $X=\Re$ we write $Lip(M)$ instead of $Lip(M,\Re)$.)

A subset $S\subset M$ will be regarded as a metric (sub-) space equipped with
the induced metric
$d|_{S\times S}$. Hence, the notations $Lip(S,X)$, and $L(f)$ for $f\in Lip(S,X)$
are clear.

A {\em simultaneous Lipschitz extension} from $S$ to $M$ is a continuous linear operator
$T: Lip(S,X)\to Lip(M,X)$ such that
$$
Tf|_{S}=f.
$$
The set of all such $T$ is denoted by $Ext(S,M;X)$ and an (optimal) {\em extension constant}
is given by
\begin{equation}\label{e2}
\lambda(S,M;X):=\inf\{ ||T||\ :\ T\in Ext(S,M;X)\}.
\end{equation}
(This becomes $\infty$, if $Ext(S,M;X)=\emptyset$.)

It is shown in [BB] that there are rather simple metric spaces 
(e.g., metric graphs with the 
vertex degrees bounded by 3) and subspaces of these
spaces for which $Ext(S,M;\Re)=\emptyset$. The results presented below
show that nevertheless there are many subspaces in a metric space for
which the extension constants (\ref{e2}) are finite.  

In what follows we will use the following definitions.

A map $\phi: (M,d)\to (M_{1},d_{1})$ is said to be {\em $C$- Lipschitz}, if its Lipschitz
constant is bounded by a constant $C$ (and simply {\em Lipschitz}, if $L(\phi)$ is bounded).

If, in addition, $\phi$ is an injection and for all $m,m'\in M$ and given $C\geq 1$
\begin{equation}\label{e3}
C^{-1}d(m,m')\leq d_{1}(\phi(m),\phi(m'))\leq Cd(m,m'),
\end{equation}
then $\phi$ is a {\em $C$-isometric embedding} (simply {\em quasi-isometric embedding}, if
(\ref{e3}) holds for some $C$).

Note that the {\em distortion} of $\phi$ (written $dst(\phi)$) satisfies
\begin{equation}\label{e4}
dst(\phi):=L(\phi)L(\phi^{-1})\leq C^{2}.
\end{equation}

Finally, $\phi:M\to M_{1}$ is a $C$-{\em isometry}, if $\phi$ is a bijection satisfying (\ref{e3}). 

Now we present the basic definitions of the paper.
\begin{D}\label{d1}
A metric space $U$ is said to be universal with respect to simultaneous Lipschitz extensions
if for an arbitrary metric space $M$ and every subspace $S$ of $M$ isometric to a
subspace of $U$
$$
\lambda(S,M;X)\leq c(U)
$$
where $c(U)$ depends only on $U$.
\end{D}
\begin{R}\label{r1'}
{\rm In fact, in all our results related to universality we will establish a much stronger
property: if, using the notation of Definition \ref{d1},
$S$ is $C$-isometric $(C\geq 1)$ to a subspace of $U$, then 
$$
\lambda(S,M;X)\leq C^{2}c(U)
$$
with $c(U)$ depending only on $U$. This clearly implies the 
universality of $U$.}
\end{R}

A deep result by Lee and Naor, see [LN, Theorem 1.6], implies universality in this sense of every
{\em doubling metric space}. Let us recall that $M$ is {\em doubling} if there is a constant $D$
such that for each $R>0$ every ball of radius $2R$ can be covered by at most $D$ balls of radius $R$. The minimal $D$ is the {\em doubling constant} of $M$ (denoted by $D(M)$).

The aforementioned theorem states that if $M_{0}$ is a doubling subspace of a metric space $M$,
then for some numerical constant $C\geq 1$
\begin{equation}\label{e5}
\sup_{X}\lambda(M_{0},M;X)\leq C\log_{2}D(M_{0}).
\end{equation}

Since every subspace $S$ of $M_{0}$ inherits the doubling property with $D(S)\leq D(M_{0})$,
inequality (\ref{e5}) implies the universality of $M_{0}$.

The main result of this paper presents a wide class of universal metric spaces which have no such hereditary property. For its formulation we require 
\begin{D}\label{d2}
A metric space is of bounded geometry with parameters $n$, $R$, $C$ if every open ball of this
space of radius $R$ admits a $C$-isometric embedding into $\Re^{n}$.
\end{D}
\begin{R}\label{r1}
{\rm Hereafter $\Re^{n}$ is regarded as the Euclidean space endowed with the standard norm
$||x||_{2}:=\{\sum_{i=1}^{n} x_{i}^{2}\}^{1/2}$, $x=(x_{1},\dots, x_{n})$.}
\end{R}
{\bf Notation.} The class of metric spaces, 
satisfying Definition \ref{d2} is denoted by  ${\cal G}_{n}(R,C)$.

Our main result is
\begin{Th}\label{te4}
Let $M:=\oplus_{i=1}^{N}M_{i}$ where every metric space $(M_{i},d_{i})$ is a (Gromov) hyperbolic space of bounded geometry. Then $M$ is universal.
\end{Th}

Here $\oplus_{i=1}^{N} M_{i}$ is a metric space with underlying set $\prod_{i=1}^{N}M_{i}$
and metric $d:=\max_{1\leq i\leq n}d_{i}$. In the sequel we also use {\em direct $p$-products}
of these spaces ($1\leq p\leq\infty$) with the same underlying set and the metric $d_{p}$ given by
\begin{equation}\label{e6}
d_{p}(m,m'):=\left(\sum_{i=1}^{N}d_{i}^{p}(m_{i},m_{i}')\right)^{1/p}.
\end{equation}

For the convenience of the reader one recalls the Rips definition of Gromov hyperbolicity (the
background material, basic properties and examples can be found in the monographs [BH] and [Gr]).
\begin{D}\label{d5}
A geodesic metric space is $\delta$-hyperbolic $(\delta\geq 0)$ if every geodesic triangle is $\delta$-slim, that is, each side of the triangle lies in the $\delta$-neighbourhood of the union of the remaining sides.

We will say that $M$ is (Gromov) hyperbolic if it is $\delta$-hyperbolic for some $\delta\geq 0$.
\end{D}

Let us also recall that
a metric space $(M,d)$ is said to be {\em geodesic}, if every pair of points can be joined by a {\em geodesic segment}. In turn, a geodesic segment joining $m$ and $m'$ is the image of a {\em geodesic}, a curve $\gamma: [0,a]\to M$ such that
$\gamma(0)=m$, $\gamma(a)=m'$ and $d(\gamma(t),\gamma(s))=|t-s|$ for all $t,s\in [0,a]$ (in particular, $a=d(m,m')$ and also is equal to the length of $\gamma$). Finally, a geodesic triangle with vertices $m_{1}\neq m_{2}\neq m_{3}$ is the union of geodesic segments with endpoints at these points.
\begin{E}\label{ex6}
{\rm (a) The Lobachevski space $\H^{n}$ is $\delta$-hyperbolic with $\delta=\ln 3$, see, e.g., [CDP]. A straightforward computation also shows that $\H^{n}$ is of bounded geometry and belongs to ${\cal G}_{n}(R,C)$ for every $R>0$ and $C=C(n,R)$.\footnote{In the sequel $C,C_{1}$ etc denote
constants; we write $C=C(\alpha,\beta,\dots)$ if the $C$ depends {\em only} on the parameters
$\alpha,\beta,\dots$.}}\\
{\rm (b) A simply connected $n$-dimensional Riemannian manifold with sectional curvature $\kappa$
satisfying $-b^{2}\leq\kappa\leq -a^{2}<0$ for some $a,b>0$, is hyperbolic (a consequence of (a)
and Toponogov's comparison theorem) and belongs to ${\cal G}_{n}(R,C)$ with $C=C(n,R)$ for
every $R>0$ (a consequence of Rauch's comparison theorem).}\\
{\rm (c) A metric tree is 0-hyperbolic, since there are no nondegenerate triangles (cycles) in this space. If the degrees of the vertices of the tree are uniformly bounded, it belongs to ${\cal G}_{2}(R,C)$, $C=C(R)$, for any $R>0$.}\\
{\rm (d) A strongly pseudoconvex domain in $\Co^{n}$ with the Kobayashi metric is Gromov hyperbolic and of bounded geometry. More generally this holds for uniform domains in $\Re^{n}$ with the quasi-hyperbolic metric, see [BHK, Chapter 1].}
\end{E}
\begin{R}\label{r6'}
{\rm Theorem \ref{te4} is of interest only in the case of unbounded geodesic spaces $M_{i}$. In fact, a bounded geodesic space is clearly Gromov hyperbolic. If, in addition, it is of bounded
geometry, then by Lemma \ref{l2.2} below this space is doubling and its universality follows from the Lee-Naor theorem.}
\end{R}

Combining Theorem \ref{te4} with the above mentioned result of Lee and Naor [LN, Theorem 1.6] one obtains the following
\begin{C}\label{new}
Let $M:=\oplus_{i=1}^{N}M_{i}$ where every $(M_{i},d_{i})$ is either a 
doubling metric space or a Gromov hyperbolic space of bounded geometry. 
Then $M$ is universal.
\end{C} 

The proof of Theorem \ref{te4} is based on several recent results on Lipschitz embeddings and extensions and a new theorem of our own 
that will be formulated now. For this goal we need several notions.
\begin{D}\label{d8}
A Borel measure on a metric space $(M,d)$ is said to be doubling at a point $m\in M$ if every open ball centered at $m$ is of finite strictly positive measure and the doubling constant
$$
D_{m}(\mu):=\sup_{R>0}\frac{\mu(B_{2R}(m))}{\mu(B_{R}(m))}<\infty.
$$
If, in addition,
$$
D(\mu):=\sup_{m\in M}D_{m}(\mu)<\infty
$$
then $\mu$ is said to be a doubling measure.
\end{D}

Here and below $B_{R}(m)$ is the 
{\em open ball} $\{m'\in M\ :\ d(m,m')<R\}$ and $\overline{B}_{R}(m)$ is the 
{\em closed ball} $\{m'\in M\ :\ d(m,m')\leq R\}$.

A metric space endowed with a fixed doubling measure is said to be of {\em homogeneous type}; that is to say, this is a triple $(M,d,\mu)$ where $\mu$ is a doubling measure. It is known, see [CW], that
$$
\log_{2} D(M)\leq c\log_{2} D(\mu)
$$
where $c>1$ is a numerical constant.

The following definition gives a useful generalization of spaces of homogeneous type.
\begin{D}\label{d9}
A metric space $(M,d)$ with a fixed family of Borel measures $\{\mu_{m}\}_{m\in M}$ on $M$ is said to be of pointwise homogeneous type if the following holds.
\begin{itemize}
\item[(i)]
{\underline{ Uniform doubling condition}}:

For every $m\in M$, $\mu_{m}$ is doubling at $m$ and
$$
D:=\sup_{m\in M}D_{m}(\mu_{m})<\infty.
$$
\item[(ii)]
{\underline{Consistency with the metric}}:

For some constant $C>0$ and all $m_{1},m_{2}\in M$ and $R>0$
$$
|\mu_{m_{1}}-\mu_{m_{2}}|(B_{R}(m))\leq C\frac{\mu_{m}(B_{R}(m))}{R}d(m_{1},m_{2})
$$
where $m=m_{1}$ or $m_{2}$.
\end{itemize}
\end{D}

The next examples clarify this definition.
\begin{E}\label{ex10}
{\rm (a) A metric space $(M,d)$ of homogeneous type endowed with a doubling measure $\mu$ clearly satisfies Definition \ref{d9} with $C=0$ and $D=D(\mu)$.\\
(b) Let $(M,d)$ be a doubling metric space (with doubling constant $D(M)$). By the Koniagin-Vol'berg theorem [KV]
(see also [LS]) $M$ carries a doubling measure $\mu$ such that
$$
\log_{2}D(\mu)\leq c\log_{2}D(M)
$$
with some numerical constant $c\geq 1$. Hence $(M,d)$ is of homogeneous type.\\
(c) $\H^{n}$ can be equipped with a family of Borel measures satisfying
the conditions of Definition \ref{d9}, see [BSh, pp. 537-540].}
\end{E}

Finally, we need
\begin{D}\label{d11}
A family of Borel measures $\{\mu_{m}\}_{m\in M}$ on a metric space $M$ is said to be $K$-uniform $(K\geq 1)$, if for all $m_{1},m_{2}\in M$ and $R>0$
$$
\mu_{m_{1}}(B_{R}(m_{1}))\leq K\mu_{m_{2}}(B_{R}(m_{2})).
$$
\end{D}

Now, all is ready to formulate our second main result. In its formulation $(M,d_{p})$ is the direct $p$-sum $\oplus_{i=1}^{N}(M_{i},d_{i})$, see (\ref{e6}), and $(M_{i},d_{i})$ is of pointwise homogeneous type with respect to a family of Borel measures $\{\mu_{m}^{i}\}_{m\in M_{i}}$ with optimal constants $D_{i}$ and $C_{i}$, $1\leq i\leq N$, see Definition \ref{d9}.
\begin{Th}\label{te12}
Assume that $\{\mu_{m}^{i}\}_{m\in M_{i}}$ is $K_{i}$-uniform for all 
$1\leq i\leq N$. Then for the extension constant of $(M,d_{p})$, see (\ref{e2}), the following inequality holds:
\begin{equation}\label{delt}
\lambda(S,M;X)\leq c_{0}(\widetilde C_{p}+1)(\log_{2}\widetilde D+1);
\end{equation}
where 
$$
\widetilde D:=\prod_{i=1}^{N}D_{i},\ \ \ \widetilde C_{p}:=\left(\sum_{i=1}^{N}C_{i}^{q}\right)^{1/q}\prod_{i=1}^{N}K_{i},
$$
$c_{0}$ is a numerical constant and $q$ relates to $p$ by $\frac{1}{p}+\frac{1}{q}=1$.
\end{Th}

For $X=\Re$ this result is proved in the authors paper [BB, Theorem 2.25]. As an easy consequence one derives from there a special case of Theorem \ref{te12} when the
target space $X$ is complemented in its second dual $X^{**}$. In particular, the result holds for dual Banach spaces $X$
(i.e., $X=Y^{*}$ for a Banach space $Y$).
But for general $X$ the proof of Theorem 2.25 presented in [BB] needs to be modified. 
This matter will be discussed in Section 3.
\begin{R}\label{r113}
{\rm (a) It is shown in the proof that for $N=1$ a sharper inequality holds:
$$
\lambda(S,M;X)\leq c_{0}(C_{1}+1)(\log_{2}D_{1}+1)
$$
where $c_{0}$ is a numerical constant. In this case, the assumption of uniform boundedness for the 
families $\{\mu_{m}^{i}\}$ is excluded from the theorem. This assumption can be eliminated also in the case $N>1$. 
Since this requires
some additional technical consideration and enlarges substantially the
right-hand side in (\ref{delt}), we will not state this generalization of
Theorem \ref{te12}.\\
(b) It is important for some applications that the extension operator of Theorem \ref{te12} maps a Lipschitz function $f:S\to X$ into a function whose range is contained in the closure of $conv\ \! f(S)$, the convex hull of $f(S)$.
In particular, all
the above formulated results are also true when the target space is
a closed convex subset of a Banach space $X$.\\
(c) It can be seen from the proof that Theorem \ref{te12} remains true for
$M_{1}$ a space of homogeneous type (with a doubling measure $\mu$).
In this case $D_{1}=D(\mu)$, $C_{1}=0$ and we may replace
$K_{1}$ by 1 in (\ref{delt}), 
see Remarks \ref{r2.10'} and \ref{r3.12'} below.}
\end{R}
\sect{\hspace*{-1em}. Proofs of Theorem \ref{te4} and Corollary \ref{new}.}
{\bf Proof of Theorem \ref{te4}.}
We need several auxiliary results the first of which is proved in 
[NPSS, Corollary 6.2]. In the forthcoming formulations, a subset of
$M$ is said to be $\epsilon-{\em dense}$ if its distance\footnote{
The distance from $S\subset M$ to $m$ is defined by
$$
d(m,S):=\inf\{d(m,m')\ :\ m'\in S\} .
$$ } from each point of $M$ is less than $\epsilon$, and $\epsilon-{\em
separated}$ if the distance between every two distinct points of the set
is more than or equal to $\epsilon$.
\begin{Proposition}\label{l2.1}
Let $F:M\to M_{1}$ be a $C$-Lipschitz map, and $A$ be an
$\epsilon$-dense subset of $M$. Assume that there exists 
$\mu\in (0,1]$ such that for all $a,a'\in A$
$$
d_{1}(F(a),F(a'))\geq\mu C d(a,a').
$$
Assume also that $M\in {\cal G}_{n_{0}}(R_{0},C_{0})$
and that
\begin{equation}\label{eq2.1}
\mu R_{0}=64\epsilon .
\end{equation}
Then there exist an integer $N=N(n_{0},C_{0})$ and a constant 
$K=K(n_{0},C_{0},R_{0},\mu,C)$ such that $M$ admits a $K$-isometric embedding into the direct sum $M_{1}\oplus\Re^{N}$.
\ \ \ \ \ $\Box$
\end{Proposition}
\begin{Proposition}\label{l2.2}
Let $(M,d)$ be a geodesic metric space belonging to 
${\cal G}_{n_{0}}(R_{0},C_{0})$. Then for every $R>0$ there exist an 
integer $n$ and a constant $C$ such that $M\in {\cal G}_{n}(R,C)$.
\end{Proposition}
{\bf Proof.} 
We must prove that every ball $B_{R}(m)$ admits a $C$-isometric embedding into some $\Re^{n}$ where $C$ and $n$ are independent of the center $m$. To find the required embedding we choose a maximal $\epsilon$-separated set
$A_{\epsilon}$ in $M$. Due to maximality, the family of balls $B_{a}:=B_{\epsilon}(a)$, $a\in A_{\epsilon}$, covers $M$. On the other hand, the family $\widetilde B_{a}:=B_{\epsilon/2}(a)$, $a\in A_{\epsilon}$, consists of pairwise disjoint balls.
\begin{Lm}\label{le2.3}
(a) If $\epsilon\leq R_{0}/2$, then the order of the open cover $\{B_{a}\}_{a\in A_{\epsilon}}$ is at most $(4C_{0})^{n_{0}}$.\\
(b) For the same $\epsilon$ and every $a\in A_{\epsilon}$ there is a bounded linear extension operator $E_{a}: Lip(B_{a}\cap A_{\epsilon})\to Lip(B_{a})$ whose norm is bounded by $24n_{0}C_{0}^{2}$.
\end{Lm}
{\bf Proof.}
(a) By definition the order of ${\cal B}:=\{B_{a}\}_{a\in A_{\epsilon}}$ is given by
$$
ord({\cal B}):=\sup_{m\in M} card\{a\in A_{\epsilon}\ :\ m\in B_{a}\}.
$$
The union $\cup\{B_{a}\ :\ B_{a}\ni m\}$ is contained in the ball $B_{2\epsilon}(m)$. As $2\epsilon\leq R_{0}$, there is a $C_{0}$-isometric embedding $\phi$ of $B_{2\epsilon}(m)$ into the Euclidean ball $B_{\rho}(\phi(m))\subset\Re^{n_{0}}$ of radius $\rho:=2C_{0}\epsilon$. On the other hand, the family $\{\phi(\widetilde B_{a})\ :\ a\in A_{\epsilon}\}$ consists of pairwise disjoint sets. This implies that that the family
of Euclidean balls $\{B_{\rho'}(\phi(a))\}_{B_{a}\ni m}$, $\rho':=
\frac{\epsilon}{2C_{0}}$, consists of pairwise disjoint sets containing in
$B_{\rho}(\phi(m))$. Comparing the $n_{0}$-measures of the sets
$\cup\{B_{\rho'}(\phi(a))\ :\ B_{a}\ni m\}$ and $B_{\rho}(\phi(m))$ we then
get
$$
\left(\frac{\epsilon}{2C_{0}}\right)^{n_{0}}card\{a\in A_{\epsilon}\ :
\ B_{a}\ni m\}\leq (2C_{0}\epsilon)^{n_{0}}.
$$
This implies the required estimate of $ord({\cal B})$.\\
(b) Since $B_{a}\subset B_{R_{0}}(a)$, there is a $C_{0}$-isometric embedding $\phi_{a}: B_{a}\to\Re^{n_{0}}$. By the Whitney extension theorem there is a bounded linear extension operator acting from $Lip(\phi_{a}(A_{\epsilon}\cap B_{a}))$ into $Lip(\Re^{n_{0}})$ whose norm is bounded by a constant $K=K(n_{0})$; in [BB, Corollary 2.24] this constant is estimated by $24n_{0}$. Then 
compositions with $\phi_{a}^{-1}$ and $\phi_{a}$ give the required operator $E_{a}: Lip(A_{\epsilon}\cap B_{a})\to Lip(B_{a})$.
\ \ \ \ \ $\Box$

Using an appropriate Lipschitz partition of unity subordinate to the cover
$\{B_{a}\ :\ a\in A_{\epsilon}\cap B_{R}(m)\}$ of the ball $B_{R}(m)$ we paste together the operators $E_{a}$ to get a linear extension operator from $Lip(A_{\epsilon}\cap B_{R}(m))$ into $Lip(B_{R}(m))$ whose norm is bounded by a constant $k$ depending only on $ord({\cal B})$ and $\sup_{a}||E_{a}||$, see [BB, Lemma 11.3] for details. Then for the
subspaces $Lip_{0}(A_{\epsilon}\cap B_{R}(m))$ and $Lip_{0}(B_{R}(m))$
of $Lip(A_{\epsilon}\cap B_{R}(m))$ and $Lip(B_{R}(m))$
determined by the condition 
$$
f(a^{*})=0\ \ \ {\rm for\ a\ fixed}\ \ \ a^{*}\in A_{\epsilon}\cap B_{R}(m),
$$
we obtain the linear extension operator
\begin{equation}\label{eq2.1'}
E: Lip_{0}(A_{\epsilon}\cap B_{R}(m))\to Lip_{0}(B_{R}(m))\ \ \
{\rm with}\ \ \ ||E||\leq k(n_{0},C_{0}).
\end{equation}

Now we use a duality argument which requires the Banach space $K(M)$
defined as the closed linear span in $Lip(M)^{*}$ of the point evaluation functionals
$$
\delta_{M}(m)[f]:=f(m),\ \ \ m\in M.
$$
By the Kantorovich-Rubinshtein duality theorem (see, e.g., [W] and 
references therein or the Appendix in [BB])
\begin{equation}\label{eq2.1"}
K(M)^{*}=Lip_{0}(M).
\end{equation}
Also, if $S\subset M$ is a subspace containing $a^{*}$, then by the McShane extension theorem $K(S)$ is naturally identified with a closed subspace of $K(M)$ and $\delta_{M}|_{S}=\delta_{S}$. 

We apply this construction to the spaces in (\ref{eq2.1'}). Since the domain of $E$ is finite-dimensional, there exists an operator
$$
P: K(B_{R}(m))\to K(A_{\epsilon}\cap B_{R}(m))
$$
such that
$$
P^{*}=E\ \ \ {\rm and}\ \ \ ||P||=||E||\leq k(n_{0}, C_{0}).
$$
Moreover, $E$ is an extension operator and therefore $P$ is a linear {\em projection} onto  $K(A_{\epsilon}\cap B_{R}(m))$.

Next, by the McShane extension theorem
$$
||\delta_{M}(m')-\delta_{M}(m'')||_{K(M)}=d(m',m''),\ \ \
m',m''\in M.
$$
In particular, $\delta_{B_{R}(m)}$ is an isometric embedding of $B_{R}(m)$ into $K(B_{R}(m))$ and the analogous
statement holds for $\delta_{B_{R}(m)\cap A_{\epsilon}}$.

Setting now $T:=P\circ\delta_{B_{R}(m)}$ we so
define a $k(n_{0},C_{0})$-Lipschitz map of $B_{R}(m)$ into $K(A_{\epsilon}\cap B_{R}(m))$ such that for $a',a''\in A_{\epsilon}\cap B_{R}(m)$ we have
$$
||T(a')-T(a'')||_{K(A_{\epsilon}\cap B_{R}(m))}=||\delta_{B_{R}(m)}(a')-\delta_{B_{R}(m)}(a'')||_{K(B_{R}(m))}=d(a',a'').
$$

Now we are under the conditions of
Proposition \ref{l2.1} with $M:=B_{R}(m)$, $M_{1}:=K(A_{\epsilon}\cap B_{R}(m))$, $C:=k(n_{0},C_{0})$ and
$\mu:=\frac{1}{k(n_{0},C_{0})}$. Choosing here $\epsilon$ equal to  $\epsilon_{0}:=\frac{R_{0}}{64k(n_{0},C_{0})}$ we derive from this proposition the following.

There exist an integer
$N=N(n_{0},C_{0})$ and a constant $k_{1}=k_{1}(n_{0},C_{0},R_{0})$ such that
$B_{R}(m)$ admits a $k_{1}$-isometric embedding into 
$K(A_{\epsilon_{0}}\cap B_{R}(m))\oplus\Re^{N}$.

Further, note that
$$
dim\ \! K(A_{\epsilon_{0}}\cap B_{R}(m))=card(A_{\epsilon_{0}}\cap B_{R}(m))-1:=d\leq N_{1}
$$
where $N_{1}$ is independent of the choice of $m$, see, e.g., [NPSS, page 18].
Also, $K(A_{\epsilon_{0}}\cap B_{R}(m))$ is $C_{1}$-isometric to $l_{\infty}^{d}$ (considered as the
space of bounded functions on $A_{\epsilon_{0}}\cap B_{R}(m)$ equal to 0 at $a^{*}$) with $C_{1}=C_{1}(\epsilon_{0},R)$. This follows from the inequalities
$$
|f(a')-f(a'')|\leq 2||f||_{l_{\infty}^{d}}\leq
\frac{2}{\epsilon_{0}}||f||_{l_{\infty}^{d}}d(a',a''),\ \ \ a',a''\in A_{\epsilon_{0}}\cap B_{R}(m),
$$
and
$$
||f||_{l_{\infty}^{d}}:=\max_{a\in A_{\epsilon_{0}}\cap B_{R}(m)}|f(a)|\leq 2L(f)R\ .
$$
Passing to the dual spaces we get from here that $K(A_{\epsilon_{0}}\cap B_{R}(m))$ is
$C_{1}$-isometric to $l_{1}^{d}$. To finish the proof of the proposition it remains to use the natural linear quasi-isometry between $l_{1}^{d}$ and $l_{2}^{d}$ and the fact that $d\leq
N_{1}$. Together with the previous statement this implies existence of a $C$-isometric embedding of $B_{R}(m)$ into $\Re^{N+N_{1}}$ with $C=
C(n_{0},C_{0},R_{0},R)$.
\ \ \ \ \ $\Box$
\begin{Lm}\label{l2.3}
Let $(M,d)$ be the direct sum $\oplus_{i=0}^{N}M_{i}$ where 
$M_{i}=\H^{n_{i}}$ for $1\leq i\leq N$, and $M_{0}$ is an 
$n_{0}$-dimensional Banach space. Then for the extension constant of $M$ we 
have
$$
\lambda(S,M,X)\leq c(M).
$$
\end{Lm}
{\bf Proof.} 
We apply Theorem \ref{te12} to our setting. In this case the Banach
space $M_{0}$ endowed with the Lebesgue measure $\lambda$ is clearly of
homogeneous type with parameters $D_{0}=2^{n_{0}}$ and $C_{0}=0$. Moreover,
$$
\lambda(B_{R}(m))=c(n_{0})R^{n_{0}}
$$
and therefore $\lambda$ is 1-uniform in the sense of Definition \ref{d11}.
Next, it was proved in [BSh, pp. 537-540] that there exist a metric
$\rho_{i}$ on $M_{i}$ equivalent to the hyperbolic metric of $M_{i}$ and
a family of Borel measures $\{\mu_{m}^{i}\}_{m\in M}$ such that
$(M_{i},\rho_{i})$ is of pointwise homogeneous type with respect to this
family, and, moreover, $\{\mu_{m}^{i}\}_{m\in M}$ is 1-uniform
on $(M_{i},\rho_{i})$. Then the required result follows from Theorem 
\ref{te12}.\ \ \ \ \ $\Box$
\begin{Lm}\label{l2.4}
Let $(M,d)$ be a Gromov hyperbolic space of bounded geometry.
Then there are an integer $n$, a constant $K\geq 1$ and a 
finite-dimensional Euclidean space $B$ such that $M$ admits a $K$-isometric
embedding into $\H^{n}\oplus B$.
\end{Lm}
{\bf Proof.}
By the Bonk-Schramm theorem [BS] there exists a {\em rough} 
$(C,k)$-{\em similarity} $\phi$ of $M$ into some $\H^{n}$ with constants
$C\geq 1$ and $k\geq 0$. In other words, $\phi:M\to\H^{n}$ satisfies for
all $m,m'$
\begin{equation}\label{e9}
Cd(m,m')-k\leq d_{h}(\phi(m),\phi(m'))\leq Cd(m,m')+k;
\end{equation}
here $d_{h}$ is the inner metric on $\H^{n}$.

For $k=0$ this implies that $C^{-1}\phi$ is  even an isometric embedding into
$\H^{n}$. So it remains to consider the case $k>0$.

Set $\epsilon:=\frac{2k}{C}$ and define $A\subset M$ to be a maximal
$\epsilon$-separated set. That is, for all $a,a'\in A$ with
$a\neq a'$
\begin{equation}\label{e10}
d(a,a')\geq\epsilon
\end{equation}
and because of maximality for every $m\in M$ there is $a\in A$ such that $d(m,a)<\epsilon$. 
From (\ref{e9}), (\ref{e10}) and the choice of $\epsilon$ 
$$
\frac{C}{2}d(a,a')\leq d_{h}(\phi(a),\phi(a'))\leq\frac{3C}{2}d(a,a')
$$
for all $a,a'\in A$. Hence, $\phi_{A}$ is a $\frac{3C}{2}$-Lipschitz map
from $A$ into $\H^{n}$. By the Lang-Pavlovi\'{c}-Schroeder extension theorem
[LPS] $\phi|_{A}$ admits a Lipschitz extension $\widehat\phi:M\to\H^{n}$
with Lipschitz constant
$$
L(\widehat\phi)\leq\frac{3}{2}c(n)C.
$$
Moreover, at points $a,a'$ of the $\epsilon$-net $A$ this map satisfies
$$
d_{h}(\widehat\phi(a),\widehat\phi(a'))=d_{h}(\phi(a),\phi(a'))\geq
\frac{C}{2}d(a,a').
$$
Finally, being a geodesic space of bounded geometry, $(M,d)$
belongs to ${\cal G}_{N}(R,C)$ for every $R>0$ and some $N$, $C$ depending
only on $R$ and the parameters in the definition of bounded geometry for $M$,
see Proposition \ref{l2.2}. Choose 
$$
R_{0}=\frac{384 k c(n)}{C}.
$$
Then, the space $(M,d)$, the $\epsilon$-separated set $A$ and the Lipschitz map 
$\widehat\phi$ satisfy the conditions of Proposition \ref{l2.1}. By this proposition
there are a constant $K\geq 1$ and a finite-dimensional Euclidean space
$B$ such that $(M,d)$ admits a $K$-isometric embedding into $\H^{n}\oplus B$.
\ \ \ \ \ $\Box$

We are now ready to finish the proof of Theorem \ref{te4}. So, let $S$ be a 
subspace of a metric space $\widehat M$ and let $\phi:S\to\oplus_{i=1}^{N}M_{i}$ be a
$C$-isometric embedding. Here $M_{i}$ is a Gromov hyperbolic space of
bounded geometry, $1\leq i\leq N$.

We must find a linear extension operator $E:Lip(S,X)\to Lip(\widehat M,X)$ whose
norm is bounded by a constant depending {\em only} on the characteristics
of the spaces $M_{i}$ and the embedding constant $C\ (\geq L(\phi))$.

For this goal we first use Lemma \ref{l2.4} to a find a $C_{1}$-isometric
embedding $\psi$ of $\oplus_{i=1}^{N}M_{i}$ into the space
$\oplus_{i=1}^{N}\H^{n_{i}}\oplus\Re^{n_{0}}$. Note that $C_{1}$ depends only 
on the characteristics of the spaces $M_{i}$. Then the composition 
$\psi\circ\phi$ is a $CC_{1}$-isometric embedding of $S$ into
$\oplus_{i=1}^{N}\H^{n_{i}}\oplus\Re^{n_{0}}$. Set
$$
\widehat S:=Image\ \!(\psi\circ\phi)\subset
\left(\bigoplus_{i=1}^{N}\H^{n_{i}}\bigoplus\Re^{n_{0}}\right)
$$
and define the linear operator $E_{1}$ on $Lip(S,X)$ by the formula
\begin{equation}\label{e12}
E_{1}f:=f\circ\phi^{-1}\circ\psi^{-1}.
\end{equation}
Then $E_{1}:Lip(S,X)\to Lip(\widehat S,X)$ and
\begin{equation}\label{e13}
||E_{1}||\leq CC_{1}.
\end{equation}

We use now Lemma \ref{l2.3} to find a linear continuous operator 
$E_{2}:Lip(\widehat S,X)\to 
Lip(\oplus_{i=1}^{N}\H^{n_{i}}\oplus\Re^{n_{0}},X)$ such that
\begin{equation}\label{e14}
E_{2}g|_{\widehat S}=g\ \ \ {\rm for}\ \ \ g\in Lip(\widehat S,X)
\end{equation}
and, in addition,
\begin{equation}\label{e15}
||E_{2}||\leq c(\overline{n})
\end{equation}
where $\overline{n}:=(n_{0},n_{1},\dots,n_{N})$.

Finally, the coordinatewise application of the Lang-Pavlovi\'{c}-Schroeder
theorem [LPS] allows us to extend the map 
$\psi\circ\phi:S\to\oplus_{i=1}^{N}\H^{n_{i}}\oplus\Re^{n_{0}}$ to a
Lipschitz map
$\Phi:\widehat M\to\oplus_{i=1}^{N}\H^{n_{i}}\oplus\Re^{n_{0}}$ such that
\begin{equation}\label{e16}
\Phi|_{S}=\psi\circ\phi\ \ \ {\rm and}\ \ \ L(\Phi)\leq c(n)CC_{1}
\end{equation}
where $n:=\sum_{i=0}^{n}n_{i}$.

Next, define the linear operator $E_{3}$ on 
$Lip(\oplus_{i=1}^{N}\H^{n_{i}}\oplus\Re^{n_{0}},X)$ by
$$
E_{3}h:=h\circ\Phi.
$$
Then $Lip(\widehat M,X)$ is the target space of $E_{3}$ and
\begin{equation}\label{e17}
||E_{3}||\leq L(\Phi)\leq c(n)CC_{1}.
\end{equation}
Moreover, by (\ref{e16})
\begin{equation}\label{e18}
(E_{3}h)|_{S}=h(\Phi|_{S})=h\circ\psi\circ\phi.
\end{equation}

Finally, define the desired linear extension operator $E$ by
$$
E=E_{3}E_{2}E_{1}.
$$
According to (\ref{e12}), (\ref{e14}) and (\ref{e18}) $E$ acts from
$Lip(S,X)$ into $Lip(\widehat M,X)$ and
$$
Ef|_{S}=f.
$$
In addition, (\ref{e13}), (\ref{e15}) and (\ref{e17}) imply that
$$
||E||\leq C^{2}C_{1}^{2}C_{2}(\overline{n}).
$$
Hence, the extension constant $\lambda(S,\widehat M;X)$ is bounded by the constant
on the right-hand side which depends only on the characteristics of the
spaces $M_{i}$ and the embedding constant $C$ of $\phi$.\ \ \ \ \ $\Box$\\
{\bf Proof of Corollary \ref{new}.}
Let $M:=\oplus_{i=1}^{N}M_{i}$. Without loss of generality we assume that $(M_{i},d_{i})$ is doubling for $i=1$ and Gromov hyperbolic of bounded geometry for $i\geq 2$. Let $S$ be a subspace of an arbitrary metric space $\widetilde M$ and $\phi:S\to M$ be a $C$-isometric embedding. Set $M^{1}:=\oplus_{i=2}^{N}M_{i}$, so that $M=M_{1}\oplus M^{1}$. By Lemma \ref{l2.4} we embed $M^{1}$ quasi-isometrically into $H:=\oplus_{i=2}^{N}\H^{n_{i}}\oplus\Re^{n_{0}}$. Further, using the map $\delta_{M_{1}}$, see the proof of Proposition \ref{l2.2}, we embed $M_{1}$ isometrically into the predual space $K(M_{1})$ of $Lip_{0}(M_{1})$. The latter, in turn, we embed isometrically into the Banach space $l_{\infty}(B)$ where $B$ is the unit ball of $K(M_{1})$. This allows us to identify the set $M$ with its image in $l_{\infty}(B)\oplus H$ and the map $\phi:S\to M$ with a quasi-isometric embedding into this image. Then $\phi=(\phi_{1},\phi_{2})$ where $\phi_{1}:S\to l_{\infty}(B)$ and $\phi_{2}:S\to H$.

Next, by the McShane extension theorem,
$\phi_{1}$ admits a Lipschitz extension to all of $\widetilde M$ preserving its Lipschitz constant while $\phi_{2}$ can be extended to all of $\widetilde M$ with Lipschitz constant bounded by $c(\sum_{i=2}^{N}n_{i},n_{0})L(\phi_{2})$,
by the Lang-Pavlovi\'{c}-Schroeder theorem [LPS]. Hence there is a Lipschitz map $\widetilde\phi:\widetilde M\to l_{\infty}(B)\oplus H$ such that $\widetilde\phi|_{S}=\phi$ and $L(\widetilde\phi)$ is bounded by a constant $c(M)C$.

Following the arguments of the proof of Theorem \ref{te4}, we now determine 
certain bounded linear extension operators $E_{1}:Lip(\phi(S),X)\to Lip(M_{1}\oplus H,X)$ and \penalty-10000 $E_{2}: Lip(M_{1}\oplus H,X)\to Lip(l_{\infty}(B)\oplus H, X)$ with bounds of their norms depending only on the basic parameters of $M$. Setting then 
$$
E(f)[x]:=(E_{2}E_{1})(f\circ\phi^{-1})[\widetilde\phi(x)],\ \ \ x\in\widetilde M,\ \ f\in Lip(S,X),
$$ 
we obtain a linear extension operator $Lip(S,X)\to Lip(\widetilde M, X)$ whose norm is bounded by the basic parameters of $M$ and $C$. This completes the proof of the corollary.

The operator $E_{1}$ is given by Theorem \ref{te12} with $M_{1}$ being
a doubling metric space, see Remarks \ref{r2.10'} and \ref{r3.12'}.

To define $E_{2}$ we first use the Lee-Naor bounded linear extension operator $\widetilde E:Lip(M_{1},X)\to Lip(l_{\infty}(B),X)$ whose norm is controlled by the doubling constant $D(M_{1})$. Moreover, $\widetilde E$ is an averaging
operator and therefore
$$
\widetilde Ef\subset\overline{conv\ \! f(M_{1})}\ \ \ {\rm (closure\ in} \ \
X).
$$
Now for every $h\in H$ we define a linear operator $\pi_{h}:Lip(M_{1}\oplus H,X)\to Lip(M_{1},X)$ by $\pi_{h}f:=f(\cdot,h)$, and then set for $f\in Lip(M_{1}\oplus H,X)$
$$
(E_{2}f)(m,h):=(\widetilde E\pi_{h}f)(m),\ \ \ (m,h)\in l_{\infty}(B)\oplus H.
$$
By this definition
$$
\begin{array}{c}
\displaystyle
||(E_{2}f)(m_{1},h_{1})-(E_{2}f)(m_{2},h_{2})||_{X}\leq\\
\\
\displaystyle
||\widetilde E(\pi_{h_{1}}f)(m_{1})-\widetilde E(\pi_{h_{1}}f)(m_{2})||_{X}
+||\widetilde E[(\pi_{h_{1}}-\pi_{h_{2}})f](m_{2})||_{X}.
\end{array}
$$
The first term in the second line is bounded by $||\widetilde E||L(f)d_{M_{1}}(m_{1},m_{2})$ while the second one is bounded by
$$
\begin{array}{c}
\displaystyle
\sup\{||x||_{X}\ :\ x\in conv[f(\cdot, h_{1})-f(\cdot,h_{2})]\}=\\
\\
\displaystyle
\sup\{||\sum\alpha_{i}[f(m_{i},h_{1})-f(m_{i},h_{2})]||_{X}\ :\ \alpha_{i}\geq 0,\
\sum\alpha_{i}=1\ \ {\rm and}\ \ \{m_{i}\}\subset M_{1}\}.
\end{array}
$$
This supremum is clearly bounded by $L(f)d_{H}(h_{1},h_{2})$. Together with the previous this gives the required estimate of the Lipschitz constant of $E_{2}f$ in $Lip(l_{\infty}(B)\oplus H,X)$ by that of $f$.\ \ \ \ \ $\Box$
\sect{\hspace*{-1em}. Proof of Theorem \ref{te12}.}
Let $S$ be a subset of a metric space $\{M,d_{p}\}$ where
$M=\prod_{i=1}^{N}M_{i}$ and
$$
d_{p}:=\left\{\sum_{i=1}^{N}d_{i}^{p}\right\}\ \ \  (1\leq p\leq\infty).
$$
Let us recall that $(M_{i},d_{i})$ is of pointwise homogeneous type with 
respect to the family $\{\mu_{m}^{i}\}_{m\in M_{i}}$ of Borel measures on
$M_{i}$, and $D_{i}$, $C_{i}$ are, respectively, the uniform doubling 
constant and the consistency constant for this family, see Definition 
\ref{d9}. Moreover, the family $\{\mu_{m}^{i}\}_{m\in M_{i}}$ is 
$K_{i}$-uniform, see Definition \ref{d11}.

Given these we must find a linear extension operator
$E:Lip(S,X)\to Lip(M,X)$ with the required estimate of its norm.

We divide the proof into three parts. First, the required
extension operator will be constructed for a single metric space of 
pointwise homogeneous type. Then we will obtain the corresponding norm
estimate for this operator. Finally, the results obtained will be
applied to prove the required result for the direct product of the spaces
$M_{i}$, $1\leq i\leq N$.\\

{\bf A.} {\underline{\bf Extension operator}}
\\

Given a metric space $(M,d)$ of pointwise homogeneous type of Definition
\ref{d9} and a subspace $S$ we now construct an extension operator $E$ 
acting from $Lip(S,X)$ into $Lip(M,X)$ and having the desired norm 
estimate. In the construction presented below $E$ acts between {\em pointed}
Lipschitz spaces $Lip_{0}(S,X)$ and $Lip_{0}(M,X)$ that are complemented
subspaces of $Lip(S,X)$ and $Lip(M,X)$ determined by the condition
$$
f(m^{*})=0;
$$
here $m^{*}$ is a fixed point in $S$. Since there exist linear projections
on the pointed subspaces of norm one, the extension operator $E$ constructed
for these subspaces gives rise to the required linear extension operator
$\widehat E$ from $Lip(S,X)$ into $Lip(M,X)$ with $||\widehat E||=||E||$.

We prefer to work with pointed Lipschitz spaces because of the following
duality result which plays an essential role in our construction.

The space $K(M)=K(M,d)$ is defined to be the closed linear span of the
point evaluation functionals
$$
\delta_{M}(m)(f):=f(m),\ \ \ m\in M,\ \ \ f\in Lip_{0}(M),
$$
in $Lip_{0}(M)^{*}$. Then the Kantorovich-Rubinshtein duality theorem
states that
\begin{equation}\label{e3.1}
K(M)^{*}=Lip_{0}(M).
\end{equation}

The map $\delta_{M}:M\to K(M)$ is readily seen to be an isometric embedding.

Consider now the map
$$
\delta_{S}:S\to K(S).
$$
By the Dugundji extension theorem [D] there exists a continuous extension
$\widehat\delta_{S}$ of $\delta_{S}$ to the whole of $M$ satisfying
\begin{equation}\label{e3.2}
\widehat\delta_{S}(M)\subset span\ \!(\delta_{S}(S)).
\end{equation}
To apply this theorem we must assume that $S$ is closed. Clearly without
loss of generality we can accept this restriction on $S$.

Let us recall Dugundji's extension construction.

Let $\{B_{m}\}_{m\in S^{c}}$ be an open cover of the open set
$S^{c}:=M\setminus S$ by the open balls
\begin{equation}\label{e3.3}
B_{m}:=B_{r_{m}}(m),\ \ \ {\rm where}\ \ \ r_{m}:=\frac{1}{3}d(m,S).
\end{equation}
Here the distance $d(m,S)$ from a point $m$ to $S$ is defined as
$\inf_{m'\in S}d(m,m')$.

Since any metric space is paracompact, there exists a continuous partition of
unity $\{p_{\alpha}\}_{\alpha\in A}$ subordinate to the cover
$\{B_{m}\}$ whose supports $U_{\alpha}:=\{m\in S^{c}\ :\ 
p_{\alpha}(m)>0\}$ form a {\em locally finite} cover of $S^{c}$.

For every $\alpha\in A$, we now pick points
\begin{equation}\label{e3.4}
m_{1}(\alpha)\in S\ \ \ {\rm and}\ \ \ m_{2}(\alpha)\in U_{\alpha}=
supp\ \!p_{\alpha}
\end{equation}
such that
\begin{equation}\label{e3.5}
d(m_{1}(\alpha),m_{2}(\alpha))< 2d(m_{2}(\alpha),S).
\end{equation}
Such points exist, since $supp\ \!p_{\alpha}$ is contained in some ball
$B_{m}$.

The aforementioned continuous extension $\widehat\delta_{S}$ is then given
by
\begin{equation}\label{e3.6}
\widehat\delta_{S}(m):=\left\{
\begin{array}{ccc}
\delta_{S}(m)&{\rm if}&m\in S\\
\\
\displaystyle
\sum_{\alpha\in A}\delta_{S}(m_{1}(\alpha))p_{\alpha}(m)&{\rm if}
&m\in S^{c}.
\end{array}
\right.
\end{equation}
\begin{Lm}\label{l3.1}
Let $f\in Lip_{0}(S,X)$. Then the function $\widehat f:M\to X$ given by
\begin{equation}\label{e3.7}
\widehat f(m):=\left\{
\begin{array}{ccc}
f(m)&{\rm if}&m\in S\\
\\
\displaystyle
\sum_{\alpha\in A}f(m_{1}(\alpha))p_{\alpha}(m)&{\rm if}
&m\in S^{c}
\end{array}
\right.
\end{equation}
satisfies for all $m,m'\in M$ the inequality
\begin{equation}\label{e3.8}
||\widehat f(m)-\widehat f(m')||_{X}\leq 7L(f)\{d(m,m')+d(m,S)+d(m',S)\}.
\end{equation}
\end{Lm}
{\bf Proof.}
In the case $m,m'\in S$, inequality (\ref{e3.8}) (even with constant 1) is trivial, 
since $\widehat f=f$ on $S$ and $d(m,S)=d(m',S)=0$.

Let now $m\in S$ and $m'\in S^{c}$. We denote by $V_{m}$ an open ball in the 
Banach space $X$ given by the inequality
\begin{equation}\label{e3.9}
||\widehat f(m)-x||_{X}< (5d(m,m')+2d(m',S))L(f),\ \ \ x\in X.
\end{equation}
Inequality (\ref{e3.8}) in this case, clearly follows from the inclusion
\begin{equation}\label{e3.10}
\widehat f(m')\in V_{m}.
\end{equation}
Since $\widehat f(m')$ is a convex combination of the points $f(m_{1}(\alpha))$, 
$\alpha\in A_{0}$, where the finite set $A_{0}$ is given by
$$
A_{0}:=\{\alpha\in A\ :\ m'\in supp\ \!p_{\alpha}\},
$$
see (\ref{e3.7}), inclusion (\ref{e3.10}) follows from  the condition
$$
f(m_{1}(\alpha))\in V_{m},\ \ \ \alpha\in A_{0}.
$$
This, in turn, is a consequence of the inequality
\begin{equation}\label{e3.11}
d(m_{1}(\alpha),m)< 5d(m,m')+2d(m',S)
\end{equation}
and the fact that $f\in Lip_{0}(S,X)$.

To prove (\ref{e3.11}) we choose for $\alpha\in A_{0}$ a point
$m(\alpha)\in S^{c}$ so that
$$
B_{m(\alpha)}\supset supp\ \!p_{\alpha}\ \ (\ni m_{2}(\alpha)).
$$
Then $m'\in B_{m(\alpha)}$, $m\in S$, and this and (\ref{e3.3})
imply that
$$
d(m(\alpha),S)\leq d(m(\alpha),m)\leq d(m(\alpha),m')+d(m',m)\leq
\frac{1}{3}d(m(\alpha),S)+d(m,m').
$$
Hence,
\begin{equation}\label{e3.11'}
d(m(\alpha),S)\leq d(m(\alpha),m)\leq\frac{3}{2}d(m,m').
\end{equation}
Further, $m_{2}(\alpha)\in B_{m(\alpha)}$ and therefore
$$
d(m_{2}(\alpha),m)\leq d(m_{2}(\alpha),m(\alpha))+d(m(\alpha),m)\leq\frac{1}{3}
d(m(\alpha),S)+d(m(\alpha),m).
$$
Combining this with the previous inequality we obtain
$$
d(m_{2}(\alpha),m)\leq 3 d(m,m').
$$
Finally, this, (\ref{e3.5}) and (\ref{e3.11'}) together with the inequality
$$
d(m_{1}(\alpha),m)\leq d(m_{1}(\alpha),m_{2}(\alpha))+d(m_{2}(\alpha),m)
$$
give the required inequality (\ref{e3.11}).

It remains to consider the case of $m,m'\in S^{c}$. Let, for the sake
of definiteness,
\begin{equation}\label{e3.12}
d(m',S)\leq d(m,S).
\end{equation}
Given $\epsilon>0$ we pick a point $m''\in S$ satisfying
$$
d(m',m'')\leq d(m',S)+\epsilon.
$$
We write now
$$
||\widehat f(m)-\widehat f(m')||_{X}\leq ||\widehat f(m)-\widehat f(m'')||_{X}+
||\widehat f(m'')-\widehat f(m')||_{X}.
$$ 
Since $m''\in S$, we can apply the estimate obtained in the previous part of
the proof to bound the right-hand side by
$$
L(f)\{2(d(m,S)+d(m',S))+5(d(m,m'')+d(m',m''))\}.
$$
Moreover, by the choice of $m''$
$$
d(m,m'')\leq d(m,m')+d(m',m'')\leq d(m,m')+d(m',S)+\epsilon.
$$
Therefore, the sum in the curly brackets is bounded by
$$
2(d(m,S)+d(m',S))+5d(m,m')+10d(m',S)+10\epsilon.
$$
This and (\ref{e3.12}), in turn, give the required inequality (\ref{e3.8}).

The lemma has been proved.\ \ \ \ \ $\Box$

We are now ready to define the required extension operator $E$. It is given for
$f\in Lip_{0}(S,X)$ by
\begin{equation}\label{e3.13}
(Ef)(m):=\left\{
\begin{array}{ccc}
f(m)&{\rm if}&m\in S\\
\\
I(\widehat f;m, d(m))&{\rm if}&m\in S^{c}.
\end{array}
\right.
\end{equation}
Here 
\begin{equation}\label{e3.15}
d(m):=d(m,S),
\end{equation}
$\widehat f$ is defined by (\ref{e3.7}), and for a locally continuous and
locally bounded function\footnote{i.e., continuous and bounded on every bounded
subset of $M$.} $g:M\to X$ we set
\begin{equation}\label{e3.14}
I(g;m,R):=\frac{1}{\mu_{m}(B_{R}(m))}\int_{B_{R}(m)}g\ \!d\mu_{m}.
\end{equation}

Let us recall that $\{\mu_{m}\}_{m\in M}$ is the family of Borel measures
on $M$, subject to Definition \ref{d9}.

Let us show that $E$ is well-defined, that is, that the vector function
$\widehat f$ is (strongly) continuous and bounded on every bounded
subset of $M$. 

Indeed, it is well-known (see, e.g., [GK]) that for any $f\in Lip_{0}(S,X)$ there
exists a uniquely defined linear continuous map 
$\widetilde f:K(S)\to X$ such that
$$
f=\widetilde f\circ\delta_{S}.
$$  
Then by the definitions of $\widehat f$, see (\ref{e3.7}), and 
$\widehat\delta_{S}$, see (\ref{e3.6}), we have
$$
\widehat f=\widetilde f\circ\widehat\delta_{S}.
$$
Since all the functions on the right-hand side are continuous and 
locally bounded, $\widehat f$ is continuous and locally bounded on
$M$. Therefore the integral $I((\widehat f;m,d(m))$ is finite.
\begin{R}\label{r3.2}
{\rm Our construction of the operator $E$ would be much 
simpler if we could define a Borel measurable map $\phi:S^{c}\to S$
satisfying the condition
$$
d(m,\phi(m))\leq Cd(m,S),\ \ \ m\in M,
$$
with some constant $C\geq 1$ independent of $m$ and $S$. Then
$\widehat f$ in (\ref{e3.13}) would be replaced by the
composite $f\circ\phi$ for which an inequality similar to inequality
(\ref{e3.8}) of Lemma \ref{l3.1} trivially holds (with $7$ replaced by $C$).\\
Unfortunately, such $\phi$ does not exist in general even in the simplest
case of $M=\Re^{2}$, see the corresponding counter-example
in the paper [N] by P. Novikov.}
\end{R}
 
At the next stage we must estimate the norm of the constructed
extension operator. The derivation presented below leads to an expression
which contains $\max\left(\frac{l}{l-1},D(l)\right)$ where the function
$D:[1,\infty)\to\Re_{+}$ is given by
\begin{equation}\label{e3.16}
D(l):=\sup_{m\in M}\sup_{R>0}\frac{\mu_{m}(B_{lR}(m))}{\mu_{m}(B_{R}(m))}.
\end{equation}
If $D(l)=l^{\lambda}$ for some constant $\lambda\geq 1$, then the term
$\max\left(\frac{l}{l-1},D(l)\right)$ can by minimized by choosing
$l:=1+1/\lambda$. This gives the bound $O(\lambda)=O(\log_{2}D)$ as
required.

In general, we have for $l\leq 2$ only the trivial estimate $D(l)\leq D$
which does not allow to achieve the result declared in Theorem \ref{te12}
(for $N=1$).

To overcome this obstacle we replace the original metric space $(M,d)$
by a new one for which the $D(l)$ is ``almost'' $l^{\lambda}$ for
some $\lambda>1$.
Moreover, this new space, say $(\widehat M,\widehat d)$, contains
an isometric copy of $(M,d)$. Therefore the extension constants of
these spaces, see (\ref{e2}), satisfy
\begin{equation}\label{e3.17}
\sup_{S\subset M}\lambda(S,M;X)\leq\sup_{\widehat S\subset\widehat M}
\lambda(\widehat S,\widehat M; X).
\end{equation}

As soon as an appropriate bound of the right-hand side of (\ref{e3.17})
via the doubling and the consistency constants for the $M$, has been obtained
the desired inequality for $\lambda(S,M;X)$ has been established.

We will realize this program in two steps:

In subsection B, we estimate the basic parameters of $\widehat M$ via
those of $M$.

In the next subsection, we obtain the required estimate for the norm of
extension operator (\ref{e3.13}). This gives the proof of Theorem
\ref{te12} for a single space $(N=1)$.\\

{\bf B. \underline{The basic properties of the extended metric
space}}
\\

The desired metric space $(\widehat M,\widehat d)$ has underlying set
\begin{equation}\label{e3.18}
M_{n}:=M\times\Re^{n}
\end{equation}
and metric given by
\begin{equation}\label{e3.19}
d_{n}:=d\oplus\delta_{n}
\end{equation}
where $d$ is the metric of the original space $M$ and $\delta_{n}$ is the
$l_{1}$-metric of $\Re^{n}$.
\\
The integer $n\geq 2$ will be chosen later to minimize the corresponding
estimates.

We then equip the space $(M_{n},d_{n})$ with the family of measures
${\cal F}_{n}:=\{\mu_{\widetilde m}\}_{\widetilde m\in M_{n}}$ where
\begin{equation}\label{e3.20}
\mu_{\widetilde m}:=\mu_{m}\otimes\lambda_{n},\ \ \ 
\widetilde  m=(m,x)\in M\times\Re^{n};
\end{equation}
here $\lambda_{n}$ is the Lebesgue measure on $\Re^{n}$ and
${\cal F}:=\{\mu_{m}\}_{m\in M}$ is the family of pointwise doubling
measures on $(M,d)$, see Definition \ref{d9}.

It is easy to show that the $M_{n}$ equipped with the family ${\cal F}_{n}$
is of pointwise homogeneous type but we need qualitative estimates of
its basic parameters in terms of those for $(M,d)$.

This goal will be achieved in several lemmas presented below. In their
formulations, $D_{n}$ and $C_{n}$ are the doubling and
consistency constants and $D_{n}(l)$ is the dilation function for
$(M_{n},d_{n})$. The function $D_{n}(l)$ is defined as in (\ref{e3.16})
with $\mu_{m}$ replaced by measure (\ref{e3.20}). We recall also that
$D$ and $C$ are the analogous constants for $(M,d)$.
\begin{Lm}\label{l128}
Assume that $n$ is related to the doubling constant $D$ 
by
\begin{equation}\label{e1233}
n\geq [\log_{2} D]+5\ .
\end{equation}
Then we have
$$D_{n}(1+1/n)\leq \frac{6}{5}e^{4}\ .
$$
\end{Lm}
{\bf Proof.}
Note that the open ball $B_{R}(\widetilde m)$ of $M_{n}$ is the set
$$
\{(m',y)\in M\times l_{1}^{n}\ :\ d(m',m)+||x-y||_{1}<R\}.
$$
Therefore an application of Fubini's theorem
yields
\begin{equation}\label{e1234}
\mu_{\widetilde m}(B_{R}(\widetilde m))=\gamma_{n}\int_{B_{R}(m)}
(R-d(m,m'))^{n}d\mu_{m}(m')\ ;
\end{equation}
here $B_{R}(m)$ is a ball of $M$ and $\gamma_{n}$ is the volume of the 
unit $l_{1}^{n}$-ball. 

We estimate this measure with $R$ replaced by
$$
R_{n}:=(1+1/n)R\ .
$$
Split the integral in (\ref{e1234}) into one over $B_{3R/4}(m)$ and one
over the remaining part $B_{R_{n}}(m)\setminus B_{3R/4}(m)$. Denote these
integrals by $I_{1}$ and $I_{2}$. For $I_{2}$ we get from (\ref{e1234})
$$
I_{2}\leq \gamma_{n}(R_{n}-3R/4)^{n}\int_{B_{R_{n}}(m)}d\mu_{m}(m')=
\gamma_{n}
\left(\frac{1}{4}+\frac{1}{n}\right)^{n}R^{n}\mu_{m}(B_{R_{n}}(m))\ .
$$
Using the doubling constant for ${\cal F}=\{\mu_{m}\}$
we further have
$$
\mu_{m}(B_{R_{n}}(m))\leq D\mu_{m}(B_{R_{n}/2}(m))\ .
$$
Moreover, by (\ref{e1233}), $D<2^{[\log_{2}D]+1}\leq\frac{1}{16}2^{n}$.
Combining all these inequalities we obtain
\begin{equation}\label{e1235}
I_{2}\leq\gamma_{n}\frac{1}{16}2^{-n}\left(1+\frac{4}{n}\right)^{n}R^{n}
\mu_{m}(B_{R_{n}/2}(m))\ .
\end{equation}

To estimate $I_{1}$ we present its integrand (which equals to that in
(\ref{e1234}) with $R$ replaced by $R_{n}$) in the following way.
$$
\left(1+\frac{1}{n}\right)^{n}(R-d(m,m'))^{n}\left(1+
\frac{d(m,m')}{(n+1)(R-d(m,m'))}\right)^{n}\ .
$$
Since $d(m,m')\leq 3R/4$ for $m'\in B_{3R/4}(m)$, the last factor is 
at most $\left(1+\frac{3}{n+1}\right)^{n}$. Hence, we have
$$
I_{1}\leq\gamma_{n}
\left(1+\frac{1}{n}\right)^{n}\left(1+\frac{3}{n+1}\right)^{n}
\int_{B_{3R/4}(m)}(R-d(m,m'))^{n}d\mu_{m}(m')\ .
$$
Using then (\ref{e1234}) we, finally, obtain
\begin{equation}\label{e1236}
I_{1}\leq e^{4}\mu_{\widetilde m}(B_{R}(\widetilde m))\ .
\end{equation}

To estimate $D_{n}(l)$ with $l=1+1/n$ it remains to bound the fractions
$$
\widetilde I_{k}:=\frac{I_{k}}{\mu_{\widetilde m}(B_{R}(\widetilde m))}\ ,
\ \ \ k=1,2\ .
$$
For $k=2$ we estimate the denominator from below as follows. Since 
$R_{n}<2R$, we bound $\mu_{\widetilde m}(B_{R}(\widetilde m))$
from below by
$$
\begin{array}{c}
\displaystyle
\gamma_{n}\int_{B_{R_{n}/2}(m)}(R-d(m,m'))^{n}d\mu_{m}(m')\geq
\gamma_{n} 2^{-n}
\left(1-\frac{1}{n}\right)^{n}R^{n}\int_{B_{R_{n}/2}(m)}d\mu_{m}(m')=\\
\\
\displaystyle
\gamma_{n}2^{-n}\left(1-\frac{1}{n}\right)^{n}R^{n}\mu_{m}(B_{R_{n}/2}(m))\ 
.
\end{array}
$$
Combining this with (\ref{e1235}) we get
$$
\widetilde I_{2}\leq\frac{1}{16}\left(1-\frac{1}{n}\right)^{-n}
\left(1+\frac{4}{n}\right)^{n}\ .
$$
Since $\left(1-\frac{1}{n}\right)^{-n}\leq\left(1-\frac{1}{5}\right)^{-5}$ 
as $n\geq 5$, we finally obtain
$$
\widetilde I_{2}\leq\frac{1}{5}e^{4}\ .
$$
For $\widetilde I_{1}$ using (\ref{e1236}) one immediately has
$$
\widetilde I_{1}\leq e^{4}\ .
$$
Hence, we have
$$
D_{n}(1+1/n)\leq\sup_{\widetilde m, R}(\widetilde I_{1}+\widetilde I_{2})<
\frac{6}{5}e^{4}\ .\ \ \ \ \ \Box
$$

Our next auxiliary result evaluates 
the consistency constant $C_{n}$ for family
${\cal F}_{n}=\{\mu_{\widetilde m}\}$ in terms of that for
${\cal F}:=\{\mu_{m}\}$. Recall that the latter constant is the $C$ in the 
inequality
\begin{equation}\label{e1237}
|\mu_{m_{1}}-\mu_{m_{2}}|(B_{R}(m_{i}))\leq
\frac{C\mu_{m_{i}}(B_{R}(m_{i}))}{R}d(m_{1},m_{2})
\end{equation}
where $m_{1},m_{2}$ are arbitrary points of $M$ and $R>0$, and $i=1,2$.
\begin{Lm}\label{l129}
$$
C_{n}\leq\left(1+\frac{4e}{3}\right)nC\ .
$$
\end{Lm}
{\bf Proof.}
Using Fubini's theorem, rewrite (\ref{e1234}) in the form
\begin{equation}\label{e1238}
\mu_{\widetilde m}(B_{R}(\widetilde m))=\beta_{n}\int_{0}^{R}
\mu_{m}(B_{s}(m))(R-s)^{n-1}ds\ 
\end{equation}
where $\beta_{n}$ is the volume of the unit sphere in $l_{1}^{n}$.
Then for $i=1,2$ we have
$$
|\mu_{\widetilde m_{1}}-\mu_{\widetilde m_{2}}|(B_{R}(\widetilde m_{i}))
\leq\beta_{n}\int_{0}^{R}|\mu_{m_{1}}-\mu_{m_{2}}|(B_{s}(m_{i}))\cdot
(R-s)^{n-1}ds\ .
$$
Divide now the interval of integration into subintervals $[0,R/n]$ and
$[R/n,R]$ and denote the corresponding integrals over these intervals by
$I_{1}$ and $I_{2}$. It suffices to find appropriate upper bounds for
$I_{k}$. Replacing $B_{s}(m_{i})$ in $I_{1}$ by the 
bigger ball $B_{s+R/n}(m_{i})$
and applying (\ref{e1237}) we obtain
$$
I_{1}\leq 
C\left(\beta_{n}\int_{0}^{R/n}\frac{\mu_{m_{i}}(B_{s+R/n}(m_{i}))}
{s+R/n}(R-s)^{n-1}ds\right)d(m_{1},m_{2})\ .
$$
Replacing $s$ by $t=s+R/n$ we bound the expression in the brackets by
$$
\left(\beta_{n}\int_{R/n}^{2R/n}\mu_{m_{i}}(B_{t}(m_{i})) (R-t)^{n-1}dt
\right)
\max_{R/n\leq t\leq 2R/n}\frac{(R+R/n-t)^{n-1}}{t(R-t)^{n-1}}\ .
$$
Since $[R/n,2R/n]\subset [0,R]$ and the maximum $<\frac{n}{R}
\left(1+\frac{1}{n-2}\right)^{n-1}<\frac{4e}{3}\frac{n}{R}$ for $n\geq 5$, 
this and (\ref{e1238}) yield
$$
I_{1}\leq\frac{4e}{3}Cn\frac{\mu_{\widetilde m_{i}}(B_{R}(\widetilde m_{i}))}
{R}d(m_{1},m_{2})\ .
$$
For the second term we get from (\ref{e1237})
$$
I_{2}\leq 
C\left(\beta_{n}\int_{R/n}^{R}\frac{\mu_{m_{i}}(B_{s}(m_{i}))}{s}
(R-s)^{n-1}ds\right)d(m_{1},m_{2})
$$
and by (\ref{e1238}) the term in the brackets is at most 
$\mu_{\widetilde m_{i}}(B_{R}(\widetilde m_{i}))\cdot\frac{n}{R}$.
Hence, we have
$$
I_{2}\leq Cn\frac{\mu_{\widetilde m_{i}}(B_{R}(\widetilde m_{i}))}{R}
d(m_{1},m_{2})\ .
$$
Further note that $d(m_{1},m_{2})\leq 
d_{n}(\widetilde m_{1},\widetilde m_{2})$. Hence, we obtain finally
the inequality
$$
|\mu_{\widetilde m_{1}}-\mu_{\widetilde m_{2}}|(B_{R}(\widetilde m_{i}))\leq
\left(1+\frac{4e}{3}\right)nC\frac{\mu_{\widetilde m_{i}}
(B_{R}(\widetilde m_{i}))}{R}d(\widetilde m_{1},\widetilde m_{2})
$$
whence $C_{n}\leq\left(1+\frac{4e}{3}\right)nC$.\ \ \ \ \ $\Box$
\begin{Lm}\label{l3.5}
Let $A_{n}:=\frac{6}{5}e^{4}n$ and $n\geq [\log_{2}D]+6$.
Then for all $R_{2}\geq R_{1}>0$ and
$\widetilde m\in M_{n}$
$$
\mu_{\widetilde m}(B_{R_{2}}(\widetilde m))-
\mu_{\widetilde m}(B_{R_{1}}(\widetilde m))\leq A_{n}
\frac{\mu_{\widetilde m}(B_{R_{2}}(\widetilde m))}{R_{2}}(R_{2}-R_{1}).
$$
\end{Lm}
{\bf Proof.} By definition $M_{n}=M_{n-1}\times\Re$ and 
$\mu_{\widetilde m}=\mu_{\widehat m}\otimes\lambda_{1}$ where
$\widehat m\in M_{n-1}$. Then by 
Fubini's theorem we have for $0<R_{1}\leq R_{2}$
$$
\mu_{\widetilde m}(B_{R_{2}}(\widetilde m))-\mu_{\widetilde m}(B_{R_{1}}
(\widetilde m))
=2\int_{R_{1}}^{R_{2}}\mu_{\widehat m}(B_{s}(\widehat m))ds\leq
\frac{2R_{2}\mu_{\widehat m}(B_{R_{2}}(\widehat m))}{R_{2}}(R_{2}-R_{1})\ .
$$
We claim that for arbitrary $l>1$ and $R>0$
\begin{equation}\label{e1242}
R\mu_{\widehat m}(B_{R}(\widehat m))\leq\frac{lD_{n-1}(l)}{2(l-1)}
\mu_{\widetilde m}(B_{R}(\widetilde m))\ .
\end{equation}
Together with the previous inequality this will yield
$$
\mu_{\widetilde m}(B_{R_{2}}(\widetilde m))-
\mu_{\widetilde m}(B_{R_{1}}(\widetilde m))\leq\frac{lD_{n-1}(l)}{l-1}\cdot
\frac{\mu_{\widetilde m}(B_{R_{2}}(\widetilde m))}{R_{2}}(R_{2}-R_{1})\ ,
$$
Finally choose here $l=1+\frac{1}{n-1}$ and use Lemma \ref{l128}. This will give the required inequality.

Hence, it remains to establish (\ref{e1242}). By the definition of
$D_{n-1}(l)$ we have for $l>1$
$$
\mu_{\widetilde m}(B_{lR}(\widetilde m))=2l\int_{0}^{R}
\mu_{\widehat m}(B_{ls}(\widehat m))ds\leq lD_{n-1}(l)\mu_{\widetilde m}
(B_{R}(\widetilde m))\ .
$$
On the other hand, replacing $[0,R]$ by $[l^{-1}R,R]$ we also have
$$
\mu_{\widetilde m}(B_{lR}(\widetilde m))\geq
2l\mu_{\widehat m}(B_{R}(\widehat m))
(R-l^{-1}R)= 2(l-1)R\mu_{\widehat m}(B_{R}(\widehat m))\ .
$$
Combining the last two inequalities we get (\ref{e1242}).\ \ \ \ \ $\Box$\\

{\bf C. \underline{Bound for the norm of the extension operator}}
\\

Let $\widehat E$ be the extension operator defined by (\ref{e3.13})-
(\ref{e3.15}) with $(M,d,{\cal F})$ replaced by $(M_{n},d_{n},{\cal F}_{n})$.
To formulate the basic result we set
\begin{equation}\label{e3.29}
K_{n}(l):=42(A_{n}+C_{n})D_{n}(l)(l+3)
\end{equation}
where $l$ and $n$ are related by
\begin{equation}\label{e3.30}
l=1+\frac{1}{n}.
\end{equation}
\begin{Proposition}\label{p3.6}
The following estimate
$$
||\widehat E||\leq 56A_{n}+\max\left(\frac{14(l+3)}{l-1},K_{n}(l)\right)
$$
is true.
\end{Proposition}

Before beginning the proof let us note that choosing here
$$
n:=[\log_{2} D]+6
$$
and applying Lemmas \ref{l128}-\ref{l3.5} we immediately obtain the 
inequality 
\begin{equation}\label{e3.31}
||\widehat E||\leq a_{0}(C+a_{1})(\log_{2} D+6)
\end{equation}
with some numerical constants $a_{0}$ and $a_{1}$. This clearly proves
Theorem \ref{te12} for $N=1$.\\
{\bf Proof.}
We have to show that for every $\widetilde m_{1},\widetilde m_{2}\in M_{n}$
\begin{equation}\label{e148}
||(\widehat Ef)(\widetilde m_{1})-(\widehat Ef)(\widetilde m_{2})||_{X}\leq 
K||f||_{Lip(S,X)}d_{n}(\widetilde m_{1},\widetilde m_{2})
\end{equation}
where $S\subset M_{n}$ and $K$ is the constant in the inequality of 
the proposition.\\
It suffices to consider only two cases:
\begin{itemize}
\item[(a)]
$\widetilde m_{1}\in S$\ and\ $\widetilde m_{2}\not\in S$;
\item[(b)]
$\widetilde m_{1}, \widetilde m_{2}\not\in S$.
\end{itemize}
We assume without loss of generality that
\begin{equation}\label{e149}
||f||_{Lip(S,X)}=1
\end{equation}
and simplify the computations by introducing the following notations:
\begin{equation}\label{e1410}
R_{i}:=d_{n}(\widetilde m_{i})\ ,\ \mu_{i}:=\mu_{\widetilde m_{i}}\ ,\ 
B_{ij}:=B_{R_{j}}(\widetilde m_{i})\ ,
\ v_{ij}:=\mu_{i}(B_{ij})\ ,\ \ 1\leq i,j\leq 2\ .
\end{equation}
We assume also for definiteness that
\begin{equation}\label{e1411}
0< R_{1}\leq R_{2}\ .
\end{equation}
By the triangle inequality we then have
\begin{equation}\label{e1412}
0\leq R_{2}-R_{1}\leq d_{n}(\widetilde m_{1},\widetilde m_{2})\ .
\end{equation}
Further, the quantities introduced satisfy the following inequalities:
\begin{equation}\label{e1414}
v_{i2}-v_{i1}\leq\frac{A_{n}v_{i2}}{R_{2}}(R_{2}-R_{1})\ ,
\end{equation}
\begin{equation}\label{e1416}
|\mu_{1}-\mu_{2}|(B_{ij})\leq\frac{C_{n}v_{ij}}{R_{j}}
d(\widetilde m_{1},\widetilde m_{2}).
\end{equation}

Now, from inequality (\ref{e3.8}) applied to our setting and the triangle
inequality we obtain
\begin{equation}\label{e1417}
\max\{||\widetilde f(\widetilde m)||_{X}\ :\ \widetilde m\in B_{i2}\}
\leq 28 R_{2}+7(i-1)d_{n}(\widetilde m_{1},\widetilde m_{2}).
\end{equation}
here $i=1,2$ and we set
\begin{equation}\label{e1418}
\widetilde f(\widetilde m):=\widehat f(\widetilde m)-
\widehat f(\widetilde m_{1})
\end{equation}
where $\widehat f$ is the extension of $f$ given by (\ref{e3.7}).\\

We now prove (\ref{e148}) for $\widetilde m_{1}\in S$ and 
$\widetilde m_{2}\not\in S$. We begin with the evident inequality
$$
||(\widehat Ef)(\widetilde m_{2})-(\widehat Ef)(\widetilde m_{1})||_{X}=
\frac{1}{v_{22}}\left|\left|
\int_{B_{22}}\widetilde f(\widetilde m)d\mu_{2}\right|\right|_{X}
\leq\max_{B_{22}}||\widetilde f||_{X},
$$
see (\ref{e1410}) and (\ref{e1418}).
Applying (\ref{e1417}) with $i=2$ we then bound this maximum by
$28R_{2}+7d_{n}(\widetilde m_{1},\widetilde m_{2})$. But 
$\widetilde m_{1}\in S$ and so
$$
R_{2}=d_{n}(\widetilde m_{2})\leq d_{n}(\widetilde m_{1},\widetilde m_{2});
$$
therefore (\ref{e148}) holds in this case with $K=35$.

The remaining case $\widetilde m_{1},\widetilde m_{2}\not\in S$ requires 
some additional auxiliary results. For their formulations we first write
\begin{equation}\label{e1419}
(\widehat Ef)(\widetilde m_{1})-(\widehat Ef)(\widetilde m_{2}):=D_{1}+D_{2}
\end{equation}
where
\begin{equation}\label{e1420}
\begin{array}{c}
D_{1}:=I(\widetilde f; \widetilde m_{1},R_{1})-
I(\widetilde f; \widetilde m_{1},R_{2})\\
\\
D_{2}:=I(\widetilde f; \widetilde m_{1},R_{2})-
I(\widetilde f; \widetilde m_{2},R_{2})\ ,
\end{array}
\end{equation}
see (\ref{e3.14}) and (\ref{e1418}).
\begin{Lm}\label{l143}
We have
$$
||D_{1}||_{X}\leq 56A_{n}d_{n}(\widetilde m_{1},\widetilde m_{2})\ .
$$
Recall that $A_{n}$ is the constant in Lemma \ref{l3.5}.
\end{Lm}
{\bf Proof.}
By (\ref{e1420}), (\ref{e1418}) and (\ref{e1410}),
$$
D_{1}=\frac{1}{v_{11}}\int_{B_{11}}\widetilde f d\mu_{1}-
\frac{1}{v_{12}}\int_{B_{12}}\widetilde f d\mu_{1}=
\left(\frac{1}{v_{11}}-\frac{1}{v_{12}}\right)\int_{B_{11}}
\widetilde f d\mu_{1}-\frac{1}{v_{12}}\int_{B_{12}\setminus B_{11}}
\widetilde f d\mu_{1}\ .
$$
This immediately implies that
$$
||D_{1}||_{X}\leq 2\cdot \frac{v_{12}-v_{11}}{v_{12}}\cdot
\max_{B_{12}}||\widetilde f||_{X}\ .
$$
Applying now (\ref{e1414}) and (\ref{e1412}), and then 
(\ref{e1417}) with $i=1$ we get the desired estimate.\ \ \ \ \ $\Box$

To obtain a similar estimate for $D_{2}$ we will use the following
two facts.
\begin{Lm}\label{l144}
Assume that for a given $l>1$
\begin{equation}\label{e1421}
d_{n}(\widetilde m_{1},\widetilde m_{2})\leq (l-1)R_{2}\ .
\end{equation}
Let for definiteness
\begin{equation}\label{eq1421}
v_{22}\leq v_{12}\ .
\end{equation}
Then we have
\begin{equation}\label{e1422}
\mu_{2}(B_{12}\Delta B_{22})\leq 2(A_{n}+C_{n})D_{n}(l)\frac{v_{12}}{R_{2}}
d_{n}(\widetilde m_{1},\widetilde m_{2})
\end{equation}
(here $\Delta$ denotes symmetric difference of sets).
\end{Lm}
{\bf Proof.}
Set
$$
R:=R_{2}+d_{n}(\widetilde m_{1},\widetilde m_{2})\ .
$$
Then $B_{12}\cup B_{22}\subset B_{R}(\widetilde m_{1})\cap 
B_{R}(\widetilde m_{2})$ and
\begin{equation}\label{eq1223}
\mu_{2}(B_{12}\Delta B_{22})\leq (\mu_{2}(B_{R}(\widetilde m_{1}))-
\mu_{2}(B_{12}))
+(\mu_{2}(B_{R}(\widetilde m_{2}))-\mu_{2}(B_{22}))\ .
\end{equation}
The first term on the right-hand side is at most
$$
|\mu_{2}-\mu_{1}|(B_{R}(\widetilde m_{1}))+|\mu_{2}-\mu_{1}|
(B_{R_{2}}(\widetilde m_{1}))+
(\mu_{1}(B_{R}(\widetilde m_{1}))-\mu_{1}(B_{R_{2}}(\widetilde m_{1}))\ .
$$
Estimating the first two terms by the inequality for the 
consistency constant (see Definition \ref{d9})  and the third by
Lemma \ref{l3.5} we bound this sum by
$$
C_{n}\left(\frac{\mu_{1}(B_{R}(\widetilde m_{1}))}{R}+
\frac{\mu_{1}(B_{R_{2}}(\widetilde m_{1}))}{R_{2}}\right)d_{n}(\widetilde m_{1},\widetilde m_{2})+
A_{n}\frac{\mu_{1}(B_{R}(\widetilde m_{1}))}{R}\ \! (R-R_{2})\ .
$$
Moreover, $R_{2}\leq R\leq lR_{2}$ and $R-R_{2}:=
d_{n}(\widetilde m_{1},\widetilde m_{2})$,
see (\ref{e1421}); taking into account (\ref{e3.16}) for $(M_{n},d_{n})$
and the notations (\ref{e1410}) we therefore have
$$
\mu_{2}(B_{R}(\widetilde m_{1}))-\mu_{2}(B_{12})\leq
[C_{n}(D_{n}(l)+1)+A_{n}D_{n}(l)]\frac{v_{12}}{R_{2}}\ \! 
d_{n}(\widetilde m_{1},\widetilde m_{2})\ .
$$
Similarly, by Lemma \ref{l3.5} and (\ref{eq1421})
$$
\begin{array}{c}
\displaystyle \mu_{2}(B_{R}(\widetilde m_{2}))-\mu_{2}(B_{22})\leq
A_{n}\frac{\mu_{2}(B_{R}(\widetilde m_{2}))}{R}\ \! (R-R_{2})\leq\\
\\
\displaystyle A_{n}D_{n}(l)\frac{v_{22}}{R_{2}}
d_{n}(\widetilde m_{1},\widetilde m_{2})\leq
A_{n}D_{n}(l)\frac{v_{12}}{R_{2}}d_{n}(\widetilde m_{1},\widetilde m_{2})\ .
\end{array}
$$
Combining the last two estimates with (\ref{eq1223}) we get the
result. \ \ \ \ \ $\Box$
\begin{Lm}\label{l145}
Under the assumptions of the previous lemma we have
\begin{equation}\label{e1423}
v_{12}-v_{22}\leq 3(A_{n}+C_{n})D_{n}(l)\frac{v_{12}}{R_{2}}
d_{n}(\widetilde m_{1},\widetilde m_{2})\ .
\end{equation}
\end{Lm}
{\bf Proof.} By (\ref{e1410}) the left-hand side is bounded by
$$
|\mu_{1}(B_{12})-\mu_{2}(B_{12})|+\mu_{2}(B_{12}\Delta B_{22})\ .
$$
Estimating these terms by (\ref{e1416}) and (\ref{e1422}) we
get the result. \ \ \ \ \ $\Box$

We now estimate $D_{2}$ from (\ref{e1420}) beginning with
\begin{Lm}\label{l146}
Under the conditions of Lemma \ref{l144} we have
$$
||D_{2}||_{X}\leq K_{n}(l)d_{n}(\widetilde m_{1},\widetilde m_{2})
$$
where $K_{n}(l):=42(A_{n}+C_{n})D_{n}(l)(l+3)$.
\end{Lm}
{\bf Proof.} By the definition of $D_{2}$ and our notation, see
(\ref{e1420}), (\ref{e1418}) and (\ref{e1410}),
$$
\begin{array}{c}
\displaystyle
||D_{2}||_{X}:=\left|\left|\frac{1}{v_{12}}\int_{B_{12}}\widetilde f d\mu_{1}-
\frac{1}{v_{22}}\int_{B_{22}}\widetilde f d\mu_{2}\right|\right|_{X}\leq
\frac{1}{v_{12}}\int_{B_{12}}||\widetilde f||_{X}\ \! d|\mu_{1}-\mu_{2}|\ +
\\
\\
\displaystyle
\frac{1}{v_{12}}\int_{B_{12}\Delta B_{22}}||\widetilde f||_{X}\ \! d\mu_{2}
+\left|\frac{1}{v_{12}}-\frac{1}{v_{22}}\right|\int_{B_{22}}
||\widetilde f||_{X}
\ \! d\mu_{2}:=J_{1}+J_{2}+J_{3}\ .
\end{array}
$$
By (\ref{e1416}) and (\ref{e1417}) with $i=1$ 
$$
J_{1}\leq\frac{1}{v_{12}}|\mu_{1}-\mu_{2}|(B_{12}) \sup_{B_{12}}
||\widetilde f||_{X}\leq
\frac{C_{n}}{R_{2}}d_{n}(\widetilde m_{1},\widetilde m_{2})
28R_{2}=28C_{n}d_{n}(\widetilde m_{1},\widetilde m_{2})\ .
$$
In turn, by (\ref{e1422}), (\ref{e1421}) and (\ref{e1417})
$$
\begin{array}{c}
\displaystyle J_{2}\leq\frac{1}{v_{12}}\mu_{2}(B_{12}\Delta
B_{22})\sup_{B_{12}\Delta B_{22}}||\widetilde f||_{X}\leq\\
\\
\displaystyle
\frac{2(A_{n}+C_{n})D_{n}(l)}{R_{2}}
d_{n}(\widetilde m_{1},\widetilde m_{2})
(7d_{n}(\widetilde m_{1},\widetilde m_{2})+28R_{2})\leq\\
\\
\displaystyle
14(A_{n}+C_{n})D_{n}(l)(l+3)d_{n}(\widetilde m_{1},\widetilde m_{2})\ .
\end{array}
$$
Finally, (\ref{e1423}), (\ref{e1417}) and (\ref{e1421}) yield
$$
J_{3}\leq 21(A_{n}+C_{n})D_{n}(l)(l+3)
d_{n}(\widetilde m_{1},\widetilde m_{2})\ .
$$
Combining these we get the required estimate.\ \ \ \ \ $\Box$

It remains to consider the case of 
$\widetilde m_{1},\widetilde m_{2}\in M_{n}$ satisfying
the inequality
$$
d_{n}(\widetilde m_{1},\widetilde m_{2})>(l-1)R_{2}
$$
converse to (\ref{e1421}). Now the definition (\ref{e1420}) of
$D_{2}$ and (\ref{e1417}) imply that
$$
||D_{2}||_{X}\leq 2\sup_{B_{12}\cup B_{22}}||\widetilde f||_{X}\leq
2(28R_{2}+7d_{n}(\widetilde m_{1},\widetilde m_{2}))\leq 
14\left(\frac{4}{l-1}+1\right)d_{n}(\widetilde m_{1},\widetilde m_{2})\ .
$$
Combining this with the inequalities of Lemmas \ref{l143} and
\ref{l146} and equality (\ref{e1419}) we obtain the required
estimate of the Lipschitz norm of the extension operator $\widehat E$. 
Actually, we have proved that
\begin{equation}\label{e1424}
||\widehat E||\leq  56A_{n}+ \max\left(\frac{14(l+3)}{l-1},K_{n}(l)\right)
\end{equation}
where $K_{n}(l)$ is the constant in (\ref{e3.29}). This gives the
proof of Theorem \ref{te12} for $N=1$.
\begin{R}\label{r2.10'}
{\rm Let us note that in the proof of this part of Theorem \ref{te12}
the condition of $K$-uniformity was not used.}
\end{R}

{\bf D. \underline{Proof of Theorem \ref{te12} for an arbitrary $N$}}\\
\\
{\bf (1)} Let, first, $p=\infty$.
Since the metric in 
$(M,d_{\infty})$ is given by 
$d_{\infty}:=\max_{1\leq i\leq
N}d_{i}$, where $d_{i}$ is the metric on $M_{i}$,
the ball $B_{R}(m)$ of $M$ is the product
of balls $B_{R}(m_{i})$ of $M_{i}$, $1\leq i\leq N$. Therefore for
the family of measures $\{\mu_{m}\}_{m\in M}$ given by the
tensor product
\begin{equation}\label{eq1238}
\mu_{m}:=\bigotimes_{i=1}^{N}\mu_{m_{i}}^{i}\ ,\ \ \
m=(m_{1},\dots, m_{N})\ ,
\end{equation}
we get
\begin{equation}\label{eq1239}
\mu_{m}(B_{R}(m))=\prod_{i=1}^{N}\mu_{m_{i}}^{i}(B_{R}(m_{i}))\ .
\end{equation}
Hence for the dilation function (\ref{e3.16}) of the family
$\{\mu_{m}\}_{m\in M}$ we get
\begin{equation}\label{eq1240}
D(l)=\prod_{i=1}^{N} D_{i}(l)
\end{equation}
where $D_{i}$ is the dilation function of $\{\mu_{m}^{i}\}_{m\in
M_{i}}$. In particular, $\{\mu_{m}\}_{m\in M}$ satisfies the
uniform doubling condition of Definition \ref{d9} with 
$\widetilde D:=D_{1}\cdots D_{N}$.

We check that condition (ii) of this definition (i.e., consistency
with the metric) holds for this family with constant
\begin{equation}\label{eq1241}
\widetilde C_{\infty}:=\left(\prod_{i=1}^{N}
K_{i}\right)\sum_{i=1}^{N}C_{i}\ .
\end{equation}

In fact, the identity
\begin{equation}\label{eq1241'}
\mu_{m}-\mu_{\widetilde
m}=\sum_{i=1}^{N}(\otimes_{j=1}^{i-1}\mu_{\widetilde
m_{j}}^{j})\otimes (\mu_{m_{i}}^{i}-\mu_{\widetilde
m_{i}}^{i})\otimes(\otimes_{j=i+1}^{N}\mu_{m_{j}}^{j})
\end{equation}
together with (\ref{eq1239}), the consistency with the metric for
each $M_{j}$  and
$K_{j}$-uniformity of $\{\mu_{m}^{j}\}_{m\in M_{j}}$ implies that
for $\widehat m=m $ or $\widetilde m$
$$
|\mu_{m}-\mu_{\widetilde m}|(B_{R}(\widehat
m))\leq\sum_{i=1}^{N}\left(\prod_{j\neq i}
K_{j}\right)C_{i}\frac{\mu_{\widehat m}(B_{R}(\widehat m))}{R}
d_{i}(m_{i},\widetilde
m_{i})\leq\widetilde C_{\infty}\frac{\mu_{\widehat m}
(B_{R}(\widehat m))}{R}d(m,\widetilde
m).
$$

Thus $(M,d_{\infty})$ is of
pointwise homogeneous type with respect
to the family (\ref{eq1238}) with optimal constants bounded by $\widetilde D$
and $\widetilde C_{\infty}$. Hence, by the previous part of the
theorem for $N=1$ we have the required estimate for 
$\lambda(S,M;X)$ in this case.\\
{\bf (2)} Let now $1\leq p<\infty$. In this case we cannot estimate 
the optimal constants $C$ and $D$ for the space
\begin{equation}\label{eq1253}
(M,d_{p}):=\oplus_{p}\{(M_{i},d_{i})\}_{1\leq i\leq N}
\end{equation}
directly. To overcome this difficulty we use the argument of the 
proof of the previous part of Theorem \ref{te12} (for $N=1$)
and isometrically embed this space into the space
$$
(\widehat M,\widehat d):=(M,d_{p})\oplus_{1}l_{1}^{a}
$$
with a suitable $a$. Hence, a point $\widehat m\in \widehat M$ is an
$(N+a)$-tuple
$$
\widehat m:=(m,x):=(m_{1},\dots, m_{N},x_{1},\dots, x_{a})
$$
with $m\in\prod_{i=1}^{N}M_{i}$ and $x\in\Re^{a}$. Moreover, the metric
$\widehat d$ is given by
$$
\widehat d(\widehat m,\widehat m'):=\left(\sum_{i=1}^{N}
d_{i}(m_{i},m_{i}')^{p}\right)^{1/p}+\sum_{i=1}^{a}|x_{i}-x_{i}'|.
$$
We endow $\widehat M$ with the family of measures given by the tensor product
$$
\mu_{\widehat m}:=\mu_{m}\otimes\lambda_{a} ,\ \ \ \widehat m\in\widehat M,
$$
where $\lambda_{a}$ is the Lebesgue measure on $\Re^{a}$ and
$\mu_{m}:=\otimes_{i=1}^{N}\mu_{m_{i}}^{i}$. \\
We will show that $\lambda(S,\widehat M;X)$ is bounded as required in Theorem
\ref{te12}. This immediately yields the desired estimate for 
$\lambda(S,M;X)$ and completes the proof of the theorem.

To accomplish this we need
\begin{Lm}\label{le1211}
The optimal uniform doubling constant $D$ of the family $\{\mu_{m}\}_{m\in M}$
satisfies
$$
D\leq\prod_{i=1}^{N}D_{i}.
$$
Recall that $D_{i}$ is the optimal uniform doubling constant of
$\{\mu_{m_{i}}^{i}\}_{m_{i}\in M_{i}}$.
\end{Lm}
{\bf Proof} (by induction on $N$). For the $\mu_{m}$-measure of the ball
$$
B_{2R}(m):=\{m'\in M\ :\ \sum_{i=1}^{N}d_{i}(m_{i},m_{i}')^{p}\leq (2R)^{p}\}
$$
we get by Fubini's theorem:
$$
\mu_{m}(B_{2R}(m))=\int_{d^{1}<(2R)^{p}}d\mu^{1}(m')\int_{d_{1}<
(2R)^{p}-d^{1}}d\mu_{1}(m_{1}').
$$
Here we set for simplicity:
$$
d^{1}:=\sum_{i=2}^{N}d_{i}(m_{i},m_{i}')^{p},\ \ \ d_{1}:=d(m_{1},m_{1}')^{p},
\ \ \
\mu^{1}:=\bigotimes_{i=2}^{N}\mu_{m_{i}}^{i},\ \ \ \mu_{1}:=\mu_{m_{1}}^{1}.
$$
The second integral is the $\mu_{1}$-measure of the ball $B_{2\rho}(m_{1})$
where $\rho:=\sqrt[p]{R^{p}-2^{-p}d^{1}}$ which is bounded by 
$D_{1}\mu_{1}(B_{\rho}(m_{1}))$. This and Fubini's theorem imply that
$$
\begin{array}{c}
\displaystyle
\mu_{m}(B_{2R}(m))\leq D_{1}\int_{2^{-p}d^{1}<R^{p}}d\mu^{1}(m')
\int_{d_{1}<R^{p}-2^{-p}d^{1}}d\mu_{1}(m_{1}')=\\
\\
\displaystyle
D_{1}\int_{d_{1}<R^{p}}
d\mu_{1}(m_{1}')\int_{d^{1}<(2R)^{p}-2^{p}d_{1}}d\mu^{1}(m').
\end{array}
$$
By the induction hypothesis the inner integral in the right-hand side
is bounded by
$$
\left(\prod_{i=2}^{N}D_{i}\right)\mu^{1}(B_{\sqrt[p]{R^{p}-d_{1}}}
(m_{2},\dots, m_{N}))=
\prod_{i=2}^{N}D_{i}\int_{d^{1}<R^{p}-d_{1}}d\mu^{1}(m').
$$
Combining this with the previous inequality to get the required result:
$$
\mu_{m}(B_{2R}(m))\leq\left(\prod_{i=1}^{N}D_{i}\right)\mu_{m}(B_{R}(m)).
\ \ \ \ \ \Box
$$

Using Lemma \ref{le1211} we estimate now the dilation function
$D_{a}(s)$ of the family $\{\mu_{\widehat m}\}$. Recall that for $s>1$
\begin{equation}\label{eq1254}
D_{a}(s):=\sup_{\widehat m\in\widehat M}\left\{
\frac{\mu_{\widehat m}(B_{sR}(\widehat m))}{\mu_{\widehat m}
(B_{R}(\widehat m))}\right\}
\end{equation}
To this end we simply apply to this setting Lemma \ref{l128} with $D$
replaced by $\prod_{i=1}^{N}D_{i}$ and $n$ by $a$. This yields
\begin{Lm}\label{le1212}
If $a\geq [\log_{2}\prod_{i=1}^{N}D_{i}]+5$, then
$$
D_{a}(1+1/a)\leq\frac{6}{5}e^{4}.\ \ \ \ \ \Box
$$
\end{Lm}

Now we estimate the consistency constant for the family 
$\{\mu_{\widehat m}\}_{\widehat m\in\widehat M}$, see Definition \ref{d9}.
To this goal we use (\ref{eq1241'}) for $\mu_{\widehat m}-\mu_{\widehat m'}$
and then apply Fubini's theorem to get for $\widehat m'':=\widehat m$ or
$\widehat m'$
\begin{equation}\label{eq1255}
\begin{array}{c}
\displaystyle
|\mu_{\widehat m}-\mu_{\widehat m'}|(B_{R}(\widehat m''))\leq
\\
\\
\displaystyle
\sum_{i=1}^{N}\int_{\delta_{a}<R}d\lambda_{a}\int_{d^{i}<(R-\delta_{a})^{p}}
d\mu_{i}'d\mu_{i}\int_{d_{i}<(R-\delta_{a})^{p}-d^{i}}
d|\mu_{m_{i}}^{i}-\mu_{m_{i}'}^{i}| .
\end{array}
\end{equation}
Here we use the notation:
$$
\begin{array}{c}
\displaystyle
\delta_{a}:=\sum_{j=1}^{a}|x_{j}-x_{j}''|,\ \ \
d^{i}:=\sum_{j\neq i}d_{j}(m_{j}'',m_{j})^{p},\ \ \ 
d_{i}:=d(m_{i}'',m_{i})^{p},\\
\\
\displaystyle
\mu_{i}':=\bigotimes_{j<i}\mu_{m_{j}'}^{j},\ \ \ \mu_{i}:=\bigotimes_{j>i}
\mu_{m_{j}}^{j}.
\end{array}
$$
Recall that $\widehat m=(m,x)\in M\times\Re^{a}$.

The inner integral in the $i$-th term of the right-hand side of 
(\ref{eq1255}) equals \penalty-10000 $|\mu_{m_{i}}^{i}-\mu_{m_{i}'}^{i}|
(B_{\rho}(m_{i}''))$ where $\rho:=\sqrt[p]{(R-\delta_{a})^{p}-d^{i}}$.
Replacing here $\rho$ by $\rho_{a}:=\sqrt[p]{(R_{a}-\delta_{a})^{p}-d^{i}}$
with $R_{a}:=(1+\frac{1}{a})R$ and applying the consistency inequality for
$(M_{i},d_{i})$ we then bound this inner integral by
$$
\frac{C_{i}\ \!\mu_{m_{i}''}^{i}(B_{\rho_{a}}(m_{i}''))}{\rho_{a}}
d_{i}(m_{i},m_{i}').
$$
Since $d^{i}\leq (R-\delta_{a})^{p}$, the denominator here is at least
$R_{a}-R=\frac{1}{a}R$. Therefore the inner integral is bounded by 
$$
\frac{a\ \! C_{i}\ \!
d_{i}(m_{i},m_{i}')}{R}\int_{d_{i}<(R_{a}-\delta_{a})^{p}-d^{i}}
d\mu_{m_{i}''}^{i}.
$$
Inserting this in (\ref{eq1255}) and replacing there $R$ by $R_{a}$ we
get
$$
|\mu_{\widehat m}-\mu_{\widehat m'}|(B_{R}(\widehat m''))\leq
\frac{a}{R}\sum_{i=1}^{N}C_{i}d_{i}(m_{i},m_{i}')
\int_{B_{R_{a}}(\widehat m'')}d\lambda_{a}\ \!d\mu_{i}'
\ \!d\mu_{i}\ \!d\mu_{m_{i}''}^{i}.
$$
To replace in this inequality each $\mu_{m_{j}'}^{j}$ (or
$\mu_{m_{j}}^{j}$) by $\mu_{m_{j}''}^{j}$ we now use the $K_{j}$-uniformity of
the family $\{\mu_{m_{j}}^{j}\}_{m_{j}\in M_{j}}$, see Definition
\ref{d11}. Applying this to the right-hand side of the previous inequality
and recalling definition (\ref{eq1254}) we estimate the $i$-th integral there
by
$$
\begin{array}{c}
\displaystyle
\left(\prod_{i=1}^{N}K_{i}\right)\int_{B_{R_{a}}(\widehat m'')}
d\lambda_{a}d\mu_{m''}=\left(\prod_{i=1}^{N}K_{i}\right)\mu_{\widehat m''}
(B_{R_{a}}(\widehat m''))\leq
\\
\\
\displaystyle
 D_{a}(1+1/a)\left(\prod_{i=1}^{N}K_{i}\right)
\mu_{\widehat m''}(B_{R}(\widehat m'')).
\end{array}
$$
Combining with the previous inequality we get for $\widehat m''=\widehat m$
or $\widehat m'$
$$
|\mu_{\widehat m}-\mu_{\widehat m'}|(B_{R}(\widehat m''))\leq
\frac{aD_{a}(1+1/a)}{R}\left(\prod_{i=1}^{N}K_{i}\right)\left(
\sum_{i=1}^{N}C_{i}d_{i}(m_{i},m_{i}')\right)
\mu_{\widehat m''}(B_{R}(\widehat m'')).
$$
By H\"{o}lder's inequality the sum in the brackets is at most
$$
\left(\sum_{i=1}^{N}C_{i}^{q}\right)^{1/q}\left(\sum_{i=1}^{N}
d_{i}(m_{i},m_{i}')^{p}\right)^{1/p}=:
\left(\sum_{i=1}^{N}C_{i}^{q}\right)^{1/q}d_{p}(m,m');
$$ 
here $\frac{1}{p}+\frac{1}{q}=1$. Hence the consistency constant
$C$ of the family $\{\mu_{\widehat m}\}_{\widehat m\in\widehat M}$
satisfies
\begin{equation}\label{eq1257}
C\leq a D_{a}(1+1/a)\left(\prod_{i=1}^{N}K_{i}\right)
\left(\sum_{i=1}^{N}
C_{i}^{q}\right)^{1/q}.
\end{equation}

Choose now $a:=[\log_{2}\prod_{i=1}^{N}D_{i}]+6$ and use Lemma \ref{l3.5}
for the space $(\widehat M,\widehat d)$ equipped with the family
$\{\mu_{\widehat m}\}_{\widehat m\in\widehat M}$. Then we have
$$
A_{a}\leq\frac{6}{5}e^{4}\left(\log_{2}\left(\prod_{i=1}^{N}D_{i}
\right)+6\right).
$$
Combining Lemma \ref{le1212} with (\ref{eq1257}) and the above inequality
we finally obtain the required result (cf. (\ref{e1424}))
$$
\lambda(S,\widehat M;X)\leq c_{0}(\widetilde C_{p}+1)
\left(\log_{2}\left(\prod_{i=1}^{N}D_{i}\right)+1\right)
$$
with $\widetilde C_{p}:=\left(\sum_{i=1}^{N}
C_{i}^{q}\right)^{1/q}\left(\prod_{i=1}^{N}K_{i}\right)$ and 
$\frac{1}{p}+\frac{1}{q}=1$.
\ \ \ \ \ $\Box$
\begin{R}\label{r3.12'}
{\rm It is easily checked that the proof of this part of Theorem \ref{te12}
is valid for the case of $M_{1}$ a doubling metric space. In fact,
due to Koniagin-Vol'berg's theorem [KV] this space can be endowed with a
doubling measure $\mu$ and therefore Theorem \ref{te12} holds for this case
with $N=1$, see Remark \ref{r2.10'}.
If $N\geq 2$ note that since the family of doubling measures for 
$M_{1}$ consists of a single measure $\mu$, the condition of 
$K_{1}$-uniformity is not needed in the proof. Hence, in this setting
Theorem \ref{te12} holds with $D_{1}=D(\mu)$, $C_{1}=0$ and with 1 instead
of $K_{1}$.}
\end{R}

\end{document}